\documentclass[a4paper, 11pt]{article}

\usepackage{amsmath,amssymb,amsthm}
\usepackage{dsfont}
\usepackage{graphicx}
\usepackage{subcaption}
\usepackage[british]{babel}
\usepackage[utf8]{inputenc}
\usepackage{xcolor}
\usepackage[colorlinks=true,linkcolor=blue,citecolor=red]{hyperref}
\usepackage{authblk}
\usepackage{blindtext}
\usepackage{epstopdf}
\usepackage{cite}
\usepackage{bm}
\usepackage{comment}
\usepackage{algorithm}
\usepackage{tikz}
\usetikzlibrary{calc,trees,positioning,arrows,chains,shapes.geometric,%
    decorations.pathreplacing,decorations.pathmorphing,shapes,%
    matrix,shapes.symbols}
\tikzset{
>=stealth',
  punktchain/.style={
    rectangle,  
     fill=cyan!40,
    draw=black, very thick,
    text width=12em, 
    minimum height=2em, 
    text centered, 
    on chain},
  line/.style={draw, thick, <-},
  element/.style={
    tape,
    top color=white,
    bottom color=blue!50!black!60!,
    minimum width=8em,
    draw=blue!40!black!90, very thick,
    text width=10em, 
    minimum height=2.5em, 
    text centered, 
    on chain},
  every join/.style={->, thick,shorten >=1pt},
  decoration={brace},
  tuborg/.style={decorate},
  tubnode/.style={midway, right=2pt},
}

\newtheorem{remark}{Remark}
\newtheorem{proposition}{Proposition}

\usepackage{mathtools}

\newcommand{\R}{\mathbb{R}}
\newcommand{\E}{\mathbb{E}}
\newcommand{\be}{\begin{equation}}
\newcommand{\ee}{\end{equation}}
\newcommand{\z}{\mathbf{z}}
\newcommand{\nt}{n_{\textrm{TOT}}}

\allowdisplaybreaks

\begin{document}
\title{Stochastic Galerkin particle methods for kinetic equations of plasmas with uncertainties}

\author[1]{Andrea Medaglia}
\author[2]{Lorenzo Pareschi}
\author[1]{Mattia Zanella}
\affil[1]{Department of Mathematics "F. Casorati", University of Pavia, Italy}
\affil[2]{Department of Mathematics and Computer Science, University of Ferrara, Italy}

\date{}

\maketitle

\abstract{The study of uncertainty propagation is of fundamental importance in plasma physics simulations. To this end, in the present work we propose a novel stochastic Galerkin (sG) particle {method} for collisional kinetic models of plasmas under the effect of uncertainties. This class of methods is based on a generalized polynomial chaos (gPC) expansion of the particles' position and velocity. In details, we introduce a stochastic particle approximation for the Vlasov-Poisson system with a BGK term describing plasma collisions. A careful reformulation of such dynamics is needed to perform the sG projection and to obtain the corresponding system for the gPC coefficients. We show that the sG particle method preserves the main physical properties of the problem, such as conservations and positivity of the solution, while achieving spectral accuracy for smooth solutions in the random space. Furthermore,  in the fluid limit the sG particle solver is designed to possess the asymptotic-preserving property necessary to obtain a sG particle scheme for the limiting Euler-Poisson system, thus avoiding the loss of hyperbolicity typical of conventional sG methods based on finite differences or finite volumes. We tested the schemes considering the classical Landau damping problem in the presence of both small and large initial uncertain perturbations, the two stream instability and the Sod shock tube problems under uncertainties. The results show that the proposed method is able to capture the correct behavior of the system in all test cases, even when the relaxation time scale is very small.
}
\\[+.2cm]

\tableofcontents

\section{Introduction}
The construction of numerical methods in plasma physics described by kinetic equations is a challenging problem. The main difficulties arise both from the high dimensionality of the problem and from the formation of multiscale structures that must be captured by a numerical solver^^>\cite{Dimarco2008,Campospinto2014,chacon2016,Crouseilles2004,Dimarco2015,Dimarco2014,Filbet2016,Sonnendrucker2013,Zhang2017}. 
The numerical methods for plasma physics developed in the literature can be essentially divided into two groups:  approaches based on direct discretizations of the corresponding system of partial differential equations (PDEs), like finite differences and finite volumes methods^^>\cite{Dimarco2015,Crouseilles2004,Dimarco2014,Zhang2017,rossmanith2011}, and approaches based on approximations of the underlying particle dynamics at different levels, like particle-in-cell (PIC) methods^^>\cite{Campospinto2014,crestetto2012,Filbet2016,Jingwei2020,Sonnendrucker2013,Carrillo2022}. Among the the latter, {Direct Simulation Monte Carlo (DSMC)} methods allow an efficient solution of the kinetic system in the presence of collisions described by the Landau operator^^>\cite{Dimarco2008,Dimarco2010}. To the first group, in the collisional case, belong the Fourier spectral methods specifically designed for the fast solution of the Landau integral term^^>\cite{Pareschi2000,Filbet2002,PD2014}.

The presence of uncertainty in the parameters of the physical model results in an additional computational challenge as it contributes to increasing the dimensionality of the system. Recently these kinds of problems have seen considerable effort by the scientific community and the introduction of numerous numerical techniques for quantifying the uncertainty^^>\cite{CPZ19,CZ19,despres2013,Pareschi2020,Chung2020,Xiao2021,Jingwei2021,Pareschi2021,DP2020}. Among these techniques stochastic Galerkin methods based on generalized polynomial chaos expansions have shown the ability to achieve spectral accuracy in the random space for smooth solutions^^>\cite{Jingwei2018,xiu2010,xiu2002}. However, due to the intrusive nature of such approaches most physical properties are lost, among which physical conservations and positivity of the solution. In addition, as emphasized in^^>\cite{despres2013,Dai2022}, the system of coefficients defining the sG method may lead to loss of hyperbolicity for nonlinear systems of conservation laws. This is the case, in particular, of the fluid limit of the Vlasov-Poisson system for a collisional plasma leading to the Euler-Poisson system. 

In this paper, following the methodology in^^>\cite{CPZ19,CZ19,Pareschi2020,Poelle2019} we introduce stochastic Galerkin particle methods for the Vlasov-Poisson system with a BGK collision term in presence of uncertainties. 
To this end we first introduce the particle approximation in absence of uncertainty and then, after generalizing the system to the case of uncertain inputs, we  project the corresponding particle system into the space of orthogonal stochastic polynomials. As we shall see, this requires a careful reformulation of the particle dynamics. By construction, the resulting particle-based scheme preserves the main physical properties and the positivity of the reconstructed distribution function. Additionally, the use of a BGK collision term allows to achieve the asymptotic-preserving property in the fluid regime and to obtain a sG particle scheme for the limiting Euler-Poisson system that avoids the loss of hyperbolicity. 

Previous results on sG schemes that avoid the loss of hyperbolicity have been based on suitable modifications of the gPC expansions^^>\cite{despres2013,Dai2022,Xiao2021} in a finite difference/volume setting. More precisely, in^^>\cite{despres2013} an entropy closure method has been proposed, in^^>\cite{Dai2022} suitable conditions to preserve hyperbolicity have been derived, whereas in^^>\cite{Xiao2021} an hybridization strategy between sG and stochastic collocation has been used. In our approach, instead, the problem is avoided thanks to the particle structure of the solution. We mention also recent results on uncertainty quantification in plasma physics obtained numerically in^^>\cite{Chung2020} and theoretically for the Landau damping in^^>\cite{Jin2019,Jin2019b}.    

The rest of the manuscript is organized as follows. In the next section we introduce the Vlasov-Poisson system in presence of uncertainties and recall its fluid limit in the collisional case characterized by the Landau integral operator. A simplified collisional model based on the BGK operator is also introduced, which will later be used in the design of the numerical solver. Section \ref{sec:sG_particle_methods} is then devoted to the detailed description of the numerical method. First the particle scheme in absence of uncertainties is presented and then it is generalized to the case of uncertain inputs through a stochastic Galerkin projection of the particle dynamics. In Section \ref{sec:test} several numerical results are reported, showing the ability of the method to achieve spectral accuracy for smooth solutions in the random space and to correctly capture the behavior of the solution in a series of prototypical test cases. Finally the paper ends with some concluding remarks and discussing future research directions.

\section{Kinetic description of plasmas with uncertainties}
\label{sec:model}

We consider the evolution of the plasma electrons at the kinetic level described by^^>\cite{Chen1974}
\be \label{eq:FP}
\dfrac{\partial f(x,v,t,\z)}{\partial t} + v \cdot \nabla_x f(x,v,t,\z) + E(x,t,\z) \cdot \nabla_v f(x,v,t,\z) = \frac1{\epsilon}Q(f,f)(x,v,t,\z),
\ee
where the one-particle distribution function $f(x,v,t,\z)$ depends on position $x\in \R^3$, velocity $v\in \R^3$ and time $t>0$, {and $\epsilon>0$ is the Knudsen number}. The random vector $\z=(z_1,\dots,z_d)\in\R^{d_z}$ characterizes the uncertainties of the system due, e.g., to missing information on the initial states, measurements of model parameters and boundary conditions. In the following we will suppose to know the distribution of $\z$, given by $p(\z)$, i.e.
\[
\mathrm{P}(\z\in \Omega) = \int_\Omega p(\z) d\z,
\]
for any $\Omega\in\R^{d_\z}$. In equation \eqref{eq:FP}, $E(x,t,\z)$ is the self-consistent electric field given by
\be 
E(x,t,\z)=-\nabla_x \phi(x,t,\z),
\ee
where $\phi(x,t,\z)$ is the potential, solution to the Poisson equation
\be\label{eq:poisson}
\Delta_x \phi(x,t,\z) = 1-\int_{\mathbb R^3}f(x,v,t,\z)dv
\ee
subject to suitable boundary conditions. In \eqref{eq:poisson} we have considered a uniform background of motionless positive ions. 

The bilinear operator $Q(f,f)$ is the collisional operator describing the interactions between the charged particles. A classical choice is given by the Landau collision operator^^>\cite{landau1965}
\[
Q(f,f)(x,v,t,\z) = \nabla_v \cdot \int_{\mathbb R^3} \Phi(v-v_*)\left[ \nabla_v f(v,\z) f(v_*,\z) - \nabla_{v_*}f(v_*,\z)f(v,\z)\right]dv_*, 
\] 
where we omitted explicit dependence from $(x,t)$ for brevity, and $\Phi$ is a $3 \times 3$ nonnegative symmetric matrix that defines the interactions between particles. Its usual form is given by 
\[
\Phi(v) = |v|^{\gamma+2}S(v), \qquad S(v) = I - \dfrac{v \otimes v}{|v|^2}, 
\]
where $\gamma>0$ in the case of hard potentials, $\gamma=0$ in the Maxwellian case, and $\gamma<0$ in the case of soft potential. In particular, the choice $\gamma = -3$ is the so-called Coulombian case and plays a central role in plasma physics. The kernel of the collision operator is characterized by local Maxwellian distibutions
\begin{equation}
\mathcal{M}_{\rho,U,T}(x,v,\z) = \rho(x,\z) \left( \dfrac{1}{2\pi T(x,\z)} \right)^{\frac{3}{2}} \exp\left( - \dfrac{ (v - U(x,\z))^2}{2 T(x,\z)}\right),
\end{equation}
where we introduced the mass, the momentum and temperature as
\[
\rho(x,\z) =  \int_{\R^3} f(x,v,t,\z) dv,
\]
\[
\rho(x,\z)U(x,\z) = \int_{\R^3}  v f(x,v,t,\z) dv, 
\]
\[
T(x,\z)=\dfrac{1}{3\rho(x,\z)}\int_{\R^3}  |v-U(x,\z)|^2 f(x,v,t,\z) dv.
\]
In the collisionless case $\epsilon\to +\infty$ from \eqref{eq:FP} we recover the Vlasov-Poisson model, whereas passing formally to the limit $\epsilon\to 0$ we obtain $f=\mathcal{M}_{\rho,U,T}$ and 
thus, defining 
\[
W(x,\z) = \rho(x,\z) \left( \dfrac{ |U(x,\z)|^2 }{ 2 } + \dfrac{3 T(x,\z)}{2} \right),\quad p(x,\z) = \rho(x,\z)T(x,\z), 
\]
we recover the Euler-Poisson system
\be \label{eq:euler_poisson}
\begin{split}
	\partial_t \rho + \nabla_x \cdot \left( \rho U \right) & = 0 \\
	\partial_t \left( \rho U \right) + \nabla_x \cdot \left( \rho U \otimes U \right) + \nabla_x p & = \rho \nabla_x \varphi \\
	\partial_t W + \nabla_x \cdot \left(\left(W+p \right) U\right) & = \rho U \cdot \nabla_x \varphi \\
	\Delta_x \varphi & =  \rho - 1 ,
\end{split}
\ee
where we omitted the dependence {from the uncertainty and the physical variables for the sake} of simplicity.

In this work, following usual approximations of collisional processes, we consider a BGK-type collisional operator^^>\cite{Andries2002,BGK}
\be \label{eq:BGK}
Q(f,\mathcal{M})(x,v,t,\z) = \mu(\mathcal{M}_{\rho,U,T}(x,v,\z) - f(x,v,t,\z)),
\ee
with $\mu>0$ the collision frequency. Recent approaches to plasma dynamics involving BGK-type operators have been proposed in^^>\cite{Xiao2021, crestetto2012, liu2017}.

For the numerical solution of \eqref{eq:FP} with uncertainties several approaches have been proposed in the literature^^>\cite{Jingwei2018,Chung2020,Xiao2021}. In most cases, however, when an intrusive method like stochastic-Galerkin is adopted, most physical properties are lost and a particular care is required to avoid the loss of hyperbolicity of the solution in fluid regimes^^>\cite{despres2013,Dai2022,Xiao2021}. In the sequel, we will consider a different approach to the problem based on a stochastic-Galerkin projection of the corresponding particle approximation.

%


\section{The stochastic Galerkin particle method} 
\label{sec:sG_particle_methods}
In this section, we first recall particle methods for the numerical approximation of equation \eqref{eq:FP} in absence of uncertainties. Without attempting to review the very huge literature on particle methods for kinetic equations, we mention Direct Simulation Monte Carlo (DSMC)^^>\cite{Pareschi2001,Pareschi2013,pareschi2001_timerelaxed,pareschi2005} and particle-in-cell (PIC) methods for collisionless plasmas, see e.g.^^>\cite{chacon2016, Filbet2016} and the references therein. It is worth to mention that extensive efforts have been devoted also to deterministic methods, see e.g.^^>\cite{Dimarco2015,Dimarco2010,Dimarco2014,Dimarco2008,Crouseilles2004, Filbet2001, Zhang2017, Sonnendrucker2013, filbet2003} and the references therein.

Furthermore, we focus on fundamentals of  stochastic Galerkin methods and we derive the corresponding sG particle approach. We highlight that these methods have been previously studied for mean-field models of collective phenomena^^>\cite{CPZ19,CZ19,medaglia2022}, and subsequently extended to the space-homogeneous Boltzmann equation^^>\cite{Pareschi2020}. 



\subsection{The particle based method in a deterministic setting}
To introduce the particle method, we first consider equation \eqref{eq:FP} without uncertainties. We are interested in the evolution of the density $f = f(x,v,t)$, $x,v \in \mathbb R^3$, $t \ge 0$, solution to \eqref{eq:FP} and complemented with the initial condition $f(x,v,0) = f_0(x,v)$. Let us consider a time discretization of the interval $[0,T_f]$ of size $\Delta t>0$. Denoting by $f^n(x,v)$ an approximation of $f(x,v,t^n)$, with $t^n=n\Delta t$, the starting point is a splitting between the numerical solution of the BGK collision step $f^*=\mathcal{C}_{\Delta t}(f^n)$ 
\be \label{eq:collision}
\begin{cases}
	\dfrac{\partial f^*}{\partial t} = \displaystyle\frac1{\varepsilon} Q(f^*,\mathcal{M}_{\rho,U,T}), \\
	f^*(x,v,0) = f^n(x,v),
\end{cases}
\ee
and the Vlasov transport step $f^{**}=\mathcal{T}_{\Delta t}(f^*)$ 
\be \label{eq:transport}
\begin{cases}
	\dfrac{\partial f^{**}}{\partial t} + v \cdot \nabla_x f^{**} + E(x,t) \cdot \nabla_v f^{**} = 0, \\
	f^{**}(x,v,0) = f^*(x,v,\Delta t).
\end{cases}
\ee
The solution at the time $t^{n+1}$ is therefore given by the first order splitting method
\[
f^{n+1}(x,v) =  \mathcal{T}_{\Delta t}\left( \mathcal{C}_{\Delta t}(f^n)(x,v) \right).
\]
Higher order methods can be considered, such as the second order Strang splitting scheme^^>\cite{strang1968}
\[
f^{n+1}(x,v) = \mathcal{C}_{\Delta t/2} \left( \mathcal{T}_{\Delta t}\left( \mathcal{C}_{\Delta t/2}(f^n)(x,v) \right)  \right).
\] 
It should be noted, however, that close to fluid regimes the above Strang splitting degenerates to first order accuracy. This drawback can be avoided using suitable IMEX Runge-Kutta techniques which up to date are not available for particle based solvers^^>\cite{PD2014}.

To reconstruct the distribution function at each time step several approaches are possible. Possible approaches are represented by constructing a histogram of position and velocities of the set of particles in the phase space, by the {weighted rule}, where each particle is counted in a computational cell and in neighboring cells with a fractions that is proportional to the overlapping area, or by reconstruction techniques based on a convolution of empirical distribution with a suitable mollifier, see e.g.^^>\cite{HE81}. {In particular, in the following, we adopt the histogram reconstruction}. We approximate the distribution function at each time step $n$ with a sample of $N$ particles $\{x^n_i,\,v^n_i\}_{i=1}^{N}$ 
\[
{f^n_N(x,v) = \dfrac{m}{N} \sum_{i=1}^{N} \delta( x - x^n_i  ) \otimes \delta( v - v^n_i  )},
\]
{being $\delta(\cdot)$ a Dirac delta and $m=\int_{\mathbb R^3} \rho(x)dx$ the total mass. We observe that the weight $\omega=m/N$ is shared by each particle, it does not influence the dynamics but it should be taken into account in the reconstruction}. The particles are then advanced in time by solving the equations of motion derived from \eqref{eq:collision} and \eqref{eq:transport}. In the following, we describe the BGK collisional step and the transport step separately. 

\paragraph{Collision step.} \label{sec:MCBGK}
In order to present the algorithm for the BGK collision step, we rewrite \eqref{eq:collision} at discrete time steps $\Delta t$ as
\[
f^*(x,v) = \left( 1 - e^{-\nu\Delta t} \right) \mathcal{M}_{\rho,U,T}(x,v) + e^{-\nu\Delta t} f^n(x,v), 
\]
{where $\nu=\mu/\epsilon$}.
From a particle point of view, we can interpret the previous relation in a simple probabilistic setting, see^^>\cite{pareschi2001_timerelaxed, Pareschi2001}. With a probability $\left( 1 - e^{-\nu\Delta t} \right)$, at the time step $n$, a particle $(x^n_i,v^n_i)$ is replaced by a particle with position and velocity sampled from the local Maxwellian $\mathcal{M}_{\rho,U,T}$. Whereas, with probability $ e^{-\nu\Delta t}$ the particle maintains its state. 

The previous local quantities are computed with a piecewise constant approach. Let us consider $N_\ell$ equally spaced cells $I_\ell$ of a {one-dimensional spatial domain $I \subset \R$} 
so that $\cup_{\ell=1}^{N_\ell} I_\ell = I $ and $I_{\ell} \cap I_{k} = \varnothing$ for $\ell \ne k$. Then, at time $t^n$, mass, momentum and temperature are computed as
\be \label{eq:local_moment}
\begin{split}
\rho^n_\ell &= \dfrac{m}{N} \sum_{i=1}^{N} \chi(x^n_i \in I_\ell), \\
U^n_\ell &= \dfrac{m}{\rho_\ell N} \sum_{i=1}^{N} {v^n_i}\chi(x^n_i \in I_\ell), \\
T^n_\ell &= \dfrac{m}{\rho_\ell N} \sum_{i=1}^{N} ({v^n_i}-U^n_\ell)^2\chi(x^n_i \in I_\ell),
\end{split}
\ee
being $\chi(\cdot)$ is the indicator function. Hence, we may rewrite the BGK collisional process for every $i=1,\dots,N$ in a compact form as follows
\be \label{eq:BGK_particle}
v^{n+1}_i = \chi\left( \xi < e^{-\nu\Delta t} \right) v^n_i + \left( 1 -  \chi\left( \xi < e^{-\nu\Delta t} \right) \right) \sum_{\ell=1}^{N_\ell} \chi\left( x^n_i \in I_\ell \right) \left( U^n_\ell + {\tilde{v}_i} \sqrt{T^n_\ell}  \right) 
\ee
{being $\xi \sim \mathcal U([0,1])$ and $\tilde{v}_i\sim\mathcal{N}(0,1)$ a standard normally distributed random variable. In order to preserve the moments in each cell, we compute the samples $\{\tilde{v}_i\}_{i=1}^N$ according to^^>\cite{pareschi2005}. Let 
\[
\tilde{V} = \frac{1}{N}\sum_{i=1}^{N} \tilde{v}_i \qquad \tilde{E} = \frac{1}{2N}\sum_{i=1}^{N} \tilde{v}^2_i
\]
be mean and energy of the sampled particles. We further rescale the samples to guarantee zero mean $u=0$ and energy $e = 1/2$, that is
\be \label{eq:scaling_normal}
\frac{1}{N}\sum_{i=1}^{N} \frac{\tilde{v}_i - \lambda}{\tau} = u \qquad  \frac{1}{2N}\sum_{i=1}^{N} \left( \frac{\tilde{v}_i - \lambda}{\tau} \right)^2 = e,
\ee
where 
\[
\tau^2 = \frac{\tilde{E} - \tilde{V}/2}{e - u/2} \qquad \lambda = \tilde{V} - \tau u.
\]}
As observed in \cite{pareschi2005} the above scaling $\tilde{v}_i\to (\tilde{v}_i - \lambda)/\tau$ being a linear transformation does not alter the normal distribution of the particle samples and, {by construction, preserves momentum and energy in each cell.}

\begin{figure}[htb]
  \centering
  \begin{minipage}{.9\linewidth}
\begin{algorithm}[H]  
\footnotesize
\caption{\small{BGK collision step} }
\label{algorithm_DSMC_BGK}
	\begin{itemize}
	\item Compute the initial position and velocity of the particles $\{x^0_i,\,v^0_i\}_{i=1}^{N}$ by sampling from the initial distribution $f^0(x,v)$;
	\item compute the velocities $\{\tilde{v}_i\}_{i=1}^N$ by sampling a standard normal distribution $\mathcal{N}(0,1)$ according to \eqref{eq:scaling_normal};
	\item for $n=0$ to $\nt-1$, given  $\{x^n_i,\,v^n_i\}_{i=1}^{N}$:
	\begin{itemize}
			\item compute $\rho^n_\ell$, $U^n_\ell$, $T^n_\ell$ in each cell $I_\ell$, $\ell=1,\dots,N_\ell$, according to \eqref{eq:local_moment};
			\item for every particle $i=1,\dots,N$:
			\begin{itemize}
				\item select $\xi$ uniformly in $[0,1]$;
				\item compute the post interaction velocities $v^{n+1}_i$ according to \eqref{eq:BGK_particle};
			\end{itemize}						
	\end{itemize}
	\item end for.
	\end{itemize}
\end{algorithm}
\end{minipage}
\end{figure}


\paragraph{Vlasov transport step.} \label{sec:DSMC_transport}
The equations of motion of the particles obtained from \eqref{eq:transport} are the following coupled set of ODEs:
\be \label{eq:ODEtransport}
\dfrac{d x_i(t) }{dt} = v_i(t), \qquad \dfrac{d v_i(t) }{dt} = E(x_i,t),
\ee
that can be solved on the computational domain through the following three steps method complemented with boundary conditions. 

Let $\{x^n_i,v^n_i\}_{i=1}^N$ be at the time step $n$ and we indicate with  $E^{n+1/2}_\ell$ the electric field in the cell $I_\ell$ at the step $n+1/2$. We consider the following splitting approach 
\be \label{eq:three_step1}
x^{n+1/2}_i = x^n_i + v^{n}_i \frac{\Delta t }{ 2}, 
\ee
\be \label{eq:three_step2}
v^{n+1}_i = v^{n}_i + \Delta t \sum_{\ell=1}^{N_\ell} E^{n+1/2}_\ell  \chi(x^{n+1/2}_i \in I_\ell),
\ee
\be \label{eq:three_step3}
x^{n+1}_i = x^{n+1/2}_i + v^{n+1}_i\frac{ \Delta t}{ 2}.
\ee
The electric field in \eqref{eq:three_step2} is computed by solving the Poisson equation for the potential \eqref{eq:poisson} with a mesh based method {on a uniform grid}. To this end, we compute as in \eqref{eq:local_moment}  the local mass on the {uniform} grid at time $n+1/2$. The process can be summarized as follows. 
{
\begin{remark} \label{remark}
The PIC-type scheme \eqref{eq:three_step1}-\eqref{eq:three_step3} evolves each sample through the splitting strategy defined in  \cite{Sonnendrucker2013} to each sample, thus it conserves mass and momentum. The mass is conserved if no particle gets in or out of the domain, since the total mass $m$ is constant in time and does not affect the dynamics. The total momentum is conserved provided the Poisson solver is such that at the discrete level 
\be \label{eq:el_f}
\int_{I} E(x)\rho(x)dx = 0,
\ee
see Proposition 10 in \cite{Sonnendrucker2013}, Section 3.2.2. This condition is met on a periodic domain. Methods satisfying this condition can be found, e.g., in Section 4.2 of the same work.
\end{remark}
}

\begin{figure}[htb]
  \centering
  \begin{minipage}{.9\linewidth}
\begin{algorithm}[H]
\footnotesize
\caption{\small{Vlasov transport step}}
\label{algorithm_DSMC_transport} 
	\begin{itemize}
		\item Compute the initial position and velocity of the particles $\{x^0_i,\,v^0_i\}_{i=1}^{N}$ by sampling from initial distribution $f^0(x,v)$;
		\item for $n=0$ to $\nt-1$, given  $\{x^n_i,\,v^n_i\}_{i=1}^{N}$:
		\begin{itemize}
				\item compute the free transport on $\frac{\Delta t}{2}$ according to \eqref{eq:three_step1}, and apply the boundary conditions;
				\item compute the mass $\rho^{n+1/2}$ and then the electric field in each cell $\{E^{n+1/2}_\ell\}_{\ell=1}^{N_\ell}$ by solving the Poisson equation with a mesh based method {satisfying condition \eqref{eq:el_f}}; 
				\item update the velocities on $\Delta t$ according to \eqref{eq:three_step2}
				\item compute the free transport on $\frac{\Delta t}{2}$ according to \eqref{eq:three_step3}, and apply the boundary conditions;
		\end{itemize}
		\item end for.
	\end{itemize}
\end{algorithm}
\end{minipage}
\end{figure}

\subsection{Stochastic Galerkin particle based methods} \label{subsec:sG}
We consider stochastic Galerkin expansion of Algorithms \ref{algorithm_DSMC_BGK}-\ref{algorithm_DSMC_transport} for the numerical approximation of  \eqref{eq:FP} with uncertainties.  
We consider a sample of $N$ particles $x_i(t^n,\z)$, $v_i(t^n,\z)$, $i=1,\dots,N$ at time $t^n=n\Delta t$, and we expand them by their generalized polynomial chaos (gPC) expansion 
\be \label{eq:gPCE}
\begin{split}
x^n_i(t,\z) \approx x^{n,M}_i(t,\z) = \sum_{h=0}^M \hat{x}^n_{i,h}(t) \Psi_h(\z), \\ v^n_i(t,\z) \approx v^M_i(t,\z) = \sum_{h=0}^M \hat{v}^n_{i,h}(t) \Psi_h(\z),
\end{split}
\ee
being $\{\Psi_h(\z)\}_{h=0}^M$ a set of polynomials of degree less or equal to $M\in\mathbb{N}$, orthonormal with respect to the measure $p(\z)d\z$
\[
\int_{\Omega} \Psi_h(\z) \Psi_k(\z) p(\z) d\z = \E_{\z} [ \Psi_h(\cdot) \Psi_k(\cdot) ] = \delta_{hk},
\]
where $\Omega\in\R^d$ and $\delta_{hk}$ is the Kronecker delta. The polynomials $\{\Psi_h(\z)\}_{h=0}^M$ are chosen following the so-called Wiener--Askey scheme^^>\cite{xiu2010, xiu2002}, depending on the distribution of the parameters. In \eqref{eq:gPCE}, for fixed $h=0,\dots,M$, the terms $\hat{x}^n_{i,h}$ and $\hat{v}^n_{i,h}$ are the projections of the positions and velocities respectively in the space generated by the polynomial of degree $h\ge 0$ 
\be \label{eq:proj}
\begin{split}
	&\hat{x}^n_{i,h} = \int_{\Omega} x^n_i(\z) \Psi_h(\z) p(\z) d\z = \E_{\z} [  x^n_i(\cdot) \Psi_h(\cdot) ], \\
	&\hat{v}^n_{i,h} = \int_{\Omega} v^n_i(\z) \Psi_h(\z) p(\z) d\z = \E_{\z} [  v^n_i(\cdot) \Psi_h(\cdot) ].
\end{split}
\ee

\paragraph{Stochastic collision step.}
To define the sG particle algorithm of the BGK step with random inputs we consider a projection on the polynomial space defined by the uncertainties of the collision process in Section^^>\ref{sec:MCBGK}. Therefore, if we consider the discretization in $N_\ell$ cells as before, we compute at each time step $t^n$ uncertain mass, mean velocity and temperature as 
\be \label{eq:local_moment_z}
\begin{split}
\rho^n_\ell(\z) &= {\dfrac{m(\z)}{N} \sum_{i=1}^{N} \chi(x^n_i(\z) \in I_\ell),}                             \\
U^n_\ell (\z)   &= {\dfrac{m(\z)}{\rho^n_\ell(\z) N} \sum_{i=1}^{N} v^n_i(\z)\chi(x^n_i(\z) \in I_\ell),}                    \\
T^n_\ell(\z)    &= {\dfrac{m(\z)}{\rho^n_\ell(\z) N} \sum_{i=1}^{N} (v^n_i(\z)-U^n_\ell(\z))^2\chi(x^n_i(\z) \in I_\ell),}
\end{split}
\ee
that are approximated through the gPC particles' expansion as
\be \label{eq:local_moment_zM}
\begin{split}
\rho^{n,M}_\ell(\z) &= {\dfrac{m(\z)}{N} \sum_{i=1}^{N} \chi(x^{n,M}_i(\z) \in I_\ell),} \\
U^{n,M}_\ell (\z)   &= {\dfrac{m(\z)}{\rho^{n,M}_\ell(\z) N} \sum_{i=1}^{N} {v^{n,M}_i(\z)}\chi(x^{n,M}_i(\z) \in I_\ell),}  \\
T^{n,M}_\ell(\z)    &=  {\dfrac{m(\z)}{\rho^{n,M}_\ell(\z) N} \sum_{i=1}^{N} ({v^{n,M}_i(\z)}-U^{n,M}_\ell(\z))^2\chi(x^{n,M}_i(\z) \in I_\ell),}
\end{split}
\ee
{where $m(\z)=\int\rho(x,\z)dx$ is the uncertain total mass}. We substitute into \eqref{eq:BGK_particle} the gPC expansion $x^{n,M}_i(\z),\,v^{n,M}_i(\z)$ of the particles 
\be \label{eq:BGK_particle_z}
\begin{split}
v^{n+1,M}_i =& \chi\left( \xi < e^{-\nu\Delta t} \right) v^{n,M}_i \\
&+ \left( 1 -  \chi\left( \xi < e^{-\nu\Delta t} \right) \right) \sum_{\ell=1}^{N_\ell} \chi\left( x^{n,M}_i \in I_\ell \right) \left( U^{n,M}_\ell + \tilde{v}_i \sqrt{T^{n,M}_\ell} \right). 
\end{split}
\ee
Hence, the projection of \eqref{eq:BGK_particle_z} on the linear space, for each $h=0,\dots,M$, is given by 
\be \label{eq:sG_BGK}
\hat{v}^{n+1}_{i,h} = \chi\left( \xi < e^{-\nu\Delta t} \right) \hat{v}^{n}_{i,h} +  \left( 1 -  \chi\left( \xi < e^{-\nu\Delta t} \right) \right)  \sum_{\ell=1}^{N_\ell} \hat{W}^{n,\ell}_{i,h},
\ee
with the collision matrix
\be \label{eq:coll_matrix}
\hat{W}^{n,\ell}_{i,h} = \int_{\Omega}  \chi\left( x^{n,M}_i(\z) \in I_\ell \right) \left( U^{n,M}_\ell(\z) + \tilde{v}_i \sqrt{T^{n,M}_\ell(\z)} \right) \Psi_{h}(\z) p(\z) d\z
\ee
and where $\tilde v$ is a random variable with normal distribution as before. 

The integral in \eqref{eq:coll_matrix} may be computed through standard Gaussian quadrature. {For example, for a one-dimensional uncertainty} using $K$ nodes we have
\be
\hat{W}^{n,\ell}_{i,h} \approx \sum_{k=1}^{K} w_k \, \chi\left( x^{n,M}_i(z_k) \in I_\ell \right) \left( U^{n,M}_\ell(z_k) + \tilde{v}_i \sqrt{T^{n,M}_\ell(z_k)} \right)\Psi_{h}(z_k)
\ee
being $\{z_k,w_k\}_{k=1}^{K}$ respectively the nodes and the weights. We remark that the local quantities in the Gaussian nodes are computed {accordingly to \eqref{eq:local_moment_zM} and the shift and scale method defined by \eqref{eq:scaling_normal} permits to preserve exactly mass, momentum and energy in each cell and at each collocation node}.
\begin{figure}[htb]
  \centering
  \begin{minipage}{.9\linewidth}
\begin{algorithm}[H]
  \footnotesize
\caption{\small{sG BGK collision step}}
\label{algorithm_sG_DSMC_BGK} 
	\begin{itemize}
		\item Compute the initial gPC expansions $\{x^{M,0}_i,\,v^{M,0}_i\}_{i=1}^N$  from the initial distribution $f^0(x,v)$;
		\item compute the velocities $\{\tilde{v}_i\}_{i=1}^N$ by sampling a standard normal distribution $\mathcal{N}(0,1)$ according to \eqref{eq:scaling_normal};
		\item for $n=0$ to $\nt-1$, given  the projections $\{\hat{x}^n_{i,h},\,\hat{v}^n_{i,h},\,i=1,\dots,N,\,h=0,\dots,M\}$:
		\begin{itemize}
			\item compute $\rho^{n,M}_\ell(z_k)$, $U^{n,M}_\ell(z_k)$,  $T^{n,M}_\ell(z_k)$ in each cell $I_\ell$, with $\ell=1,\dots,N_\ell$, $k = 1,\dots,K$ according to \eqref{eq:local_moment_zM};
			\item for every particle $i=1,\dots,N$:
			\begin{itemize}
				\item select $\xi$ uniformly in $(0,1)$;
				\item compute the post interaction projections $\hat{v}^{n+1}_{i,h}$ according to \eqref{eq:sG_BGK};
			\end{itemize}						
		\end{itemize}
		\item end for.
	\end{itemize}
\end{algorithm}
\end{minipage}
\end{figure}

\paragraph{Stochastic Vlasov transport step.}
In the presence of uncertainties the equations of motion of the particles are given by 
\[
\dfrac{dx_i(t,\z)}{dt} = v_i(t,\z), \qquad  \dfrac{dv_i(t,\z)}{dt} = E(x_i,t,\z). 
\]
The gPC expansion of the particles' systems  $x^M_i(t,\z),\,v^M_i(t,\z)$ is solution to
\[
\dfrac{dx_i^M(t,\z)}{dt} = v_i^M(t,\z), \qquad  \dfrac{dv_i^M(t,\z)}{dt} = E^M(x_i^M,t,\z). 
\]
Hence, we project the latter set of ODEs in the linear space $\{\Psi_h(\z)\}_{h=0}^M$ to obtain
\be \label{eq:ODEtransport_proj}
\dfrac{d \hat{x}_{i,h}(t) }{dt} = \hat{v}_{i,h}(t), \qquad \dfrac{d \hat{v}_{i,h}(t) }{dt} = \int_{\Omega} E^M(x^{M}_i,t,\z) \Psi_h(\z)p(\z)d\z  .
\ee
Considering the time and spatial discretization introduced before, given the projections $\{\hat{x}^n_{i,h},\hat{v}^n_{i,h}\}$ and the electric field $\{E^{n+1/2,M}_\ell(\z)\}_{\ell=1}^{N_\ell}$, the three step method at the level of projected quantities reads
\be \label{eq:sG_three_step1}
\hat{x}^{n+1/2}_{i,h} = \hat{x}^n_{i,h} + \hat{v}^{n}_{i,h} \Delta t / 2, 
\ee
\be \label{eq:sG_three_step2}
\hat{v}^{n+1}_{i,h} = \hat{v}^{n}_{i,h} + \Delta t \sum_{\ell=1}^{N_\ell} \int_{\Omega} E^{n+1/2,M}_\ell(\z)  \chi(x^{n+1/2,M}_i(\z) \in I_\ell) \Psi_h(\z) p(\z) d\z,
\ee
\be \label{eq:sG_three_step3}
\hat{x}^{n+1}_{i,h} = \hat{x}^{n+1/2}_{i,h} + \hat{v}^{n+1}_{i,h} \Delta t / 2.
\ee
The integral in \eqref{eq:sG_three_step2} is again computed with Gauss quadrature on $K$ nodes: as a consequence, the electric field needs to be calculated for every Gaussian nodes. This can be done by solving in parallel the Poisson equation as exposed in Section \ref{sec:DSMC_transport} for every $\{z_k\}_{k=1}^K$.

{
Since $m(\z)$ does not depend on time and does not affect the dynamics, the scheme conserves the total mass. Moreover, the following proposition holds. 
\begin{proposition}
The sG Vlasov transport step conserves the total momentum provided the Poisson solver is such that 
\be
\sum_{\ell=1}^{N_\ell} E^{n+1/2,M}_\ell(\z)  \rho^{n+1/2,M}_{\ell}(\z) = 0,
\label{eq:assP}
\ee
at each collocation node used to evaluate \eqref{eq:sG_three_step2}.
\end{proposition}
\begin{proof}
Summing up the contributions of all the particles we obtain from  \eqref{eq:sG_three_step2} 
\[
\begin{split}
\sum_{i=1}^{N}\hat{v}^{n+1}_{i,h} &= \sum_{i=1}^{N}\hat{v}^{n}_{i,h} + \Delta t \sum_{\ell = 1}^{N_\ell} \int_{\Omega} E_\ell^{n+1/2,M}(\z)\sum_{i=1}^N \chi(x_i^{n+1/2,M}(\z)\in I_\ell) \Psi_h(\z)p(\z)d\z \\
&=\sum_{i=1}^{N}\hat{v}^{n}_{i,h}  + \Delta t \int_{\Omega} \frac{N}{m(\z)} \left( \sum_{\ell=1}^{N_\ell} E^{n+1/2,M}_\ell(\z)  \rho^{n+1/2,M}_{\ell}(\z) \right) \Psi_h(\z) p(\z) d\z,
\end{split}
\]
with $\rho_\ell^M(\z)$ defined in \eqref{eq:local_moment_zM}. Now, since by assumption \eqref{eq:assP} the term within brackets vanishes at the collocation nodes used to evaluate the above integral, we have 
\[
\sum_{i=1}^{N}\hat{v}^{n+1}_{i,h} = \sum_{i=1}^{N}\hat{v}^{n}_{i,h}
\]
 which implies moment conservation. 
 \end{proof}
}

\begin{figure}[htb]
  \centering
  \begin{minipage}{.9\linewidth}
\begin{algorithm}[H] 
\footnotesize
\caption{\small{sG Vlasov transport step}}
\label{algorithm_sG_DSMC_transport} 
	\begin{itemize}
		\item Compute the initial gPC expansions $\{x^{0,M}_i,v^{0,M}_i\}_{i=1}^N$ from the initial distribution $f^0(x,v)$;
		\item for $n=0$ to $\nt-1$, \\
		given  the projections $\{(\hat{x}^n_{i,h},\,\hat{v}^n_{i,h}),\,i=1,\dots,N,\,h=0,\dots,M\}$:
		\begin{itemize}
			\item compute the free transport on $\Delta t / 2$ according to \eqref{eq:sG_three_step1}, and apply the boundary conditions;
			\item compute the mass $\rho^{n+1/2}_l(z_k)$ in the nodes and then the electric field $\{E^{n+1/2}_l(z_k)\}_{l=1}^{N_l}$ by solving the Poisson equation with a mesh based method {satisfying condition \eqref{eq:assP}}; 
			\item update the velocities on $\Delta t$ according to \eqref{eq:sG_three_step2}
			\item compute the free transport on $\Delta t / 2$ according to \eqref{eq:sG_three_step3}, and apply the boundary conditions;
		\end{itemize}
		\item end for.
	\end{itemize}
\end{algorithm}
	\end{minipage}
\end{figure}

In the following we are interested both in periodic and reflecting boundary conditions on the space domain {$I=[L_-,L_+]$} discretized by means of $N_\ell$ cells. To this end, we may observe that periodic conditions correspond to rewrite the position of the $i$th particle at time $t^n$ as follows
{
\[
x_i^{n,M}(\z) = 
\begin{cases}
x_i^{n,M}(\z) & L_- \le x^{n,M}_i(\z) \le L_+ \\
x^{n,M}_i(\z) + (L_{+}-L_+) &   x^{n,M}_i(\z) <  L_- \\
x^{n,M}_i(\z) - (L_{+}-L_+)   & x^{n,M}_i(\z) > L_+. 
\end{cases}
\]
Therefore, at the level of projections we have for all $h = 0,\dots,M$
\be \label{eq:pbc}
\begin{split}
	\hat{x}^{n}_{i,h} = \int_{\Omega} & \Bigg\{ \chi\left( x^{n,M}_i(\z) < L_- \right) \left(x^{n,M}_i(\z) + (L_{+}-L_-) \right) + \chi\left( x^{n,M}_i(\z) > L_+\right) \left(x^{n,M}_i(\z) - (L_{+}-L_-) \right) \\
	& + \chi\left( L_-\le x^{n,M}_i(\z) \le L_+ \right) x^{n,M}_i(\z) \Bigg\} \Psi_{h}(\z) p(\z) d\z.
\end{split}
\ee
On the other hand, reflecting boundary conditions^^>\cite{Pareschi2001} on the pair $(\tilde{x}^{n,M}_i(\z),\,\tilde{v}^{n,M}_i(\z))$ reads
\[
\begin{cases}
x_i^{n,M}(\z) = \tilde{x}_i^{n,M}(\z),\\[-.2cm]
  &\textrm{if}\quad L_- \le x^{n,M}_i(\z) \le L_+ \\[-.2cm]
v^{n,M}_i(\z)=\tilde{v}^{n,M}_i(\z)\\[+.2cm]
x_i^{n,M}(\z) = L_--|\tilde{x}_i^{n,M}(\z)+L_+|\mathrm{sgn}(\tilde{v}^{n,M}_i(\z)),\\[-.2cm]
 &\textrm{if}\quad  x^{n,M}_i(\z) < L_- \\[-.2cm]
v^{n,M}_i(\z)=-\tilde{v}^{n,M}_i(\z)\\[+.2cm]
x_i^{n,M}(\z) = L_+-|\tilde{x}_i^{n,M}(\z)L_-|\mathrm{sgn}(\tilde{v}^{n,M}_i(\z)), \\[-.2cm]
 &\textrm{if}\quad  x^{n,M}_i(\z) > L_+,\\[-.2cm]
v^{n,M}_i(\z)=-\tilde{v}^{n,M}_i(\z) 
\end{cases}
\]
from which we have the projections for all $h = 0,\dots,M$
\be \label{eq:rbc_x}
\begin{split}
	\hat{x}^{n}_{i,h} = \int_{\Omega} & \Bigg\{ \chi\left( \tilde{x}^{n,M}_i(\z) < L_-\right) \left( L_--|\tilde{x}_i^{n,M}(\z)+L_+|\mathrm{sgn}(\tilde{v}^{n,M}_i(\z))  \right) \\
	& + \chi\left( \tilde{x}^{n,M}_i(\z) > L_+\right) \left( L_+-|\tilde{x}_i^{n,M}(\z)L_-|\mathrm{sgn}(\tilde{v}^{n,M}_i(\z)) \right) \\
	& + \chi\left( L_-\le \tilde{x}^{n,M}_i(\z) \le L_+\right) \tilde{x}^{n,M}_i(\z) \Bigg\} \Psi_{h}(\z) p(\z) d\z,
\end{split}
\ee
\be \label{eq:rbc_v}
\begin{split}
\hat{v}^{n}_{i,h} =& \int_{\Omega} \Bigg\{ \chi\left( \tilde{x}^{n,M}_i(\z) < L_-\right) \left(-\tilde{v}^{n,M}_i\right) + \chi\left( \tilde{x}^{n,M}_i(\z) > L_+\right) \left(-\tilde{v}^{n,M}_i\right) \\
&\qquad + \chi\left( L_-\le \tilde{x}^{n,M}_i(\z) \le L_+\right) \tilde{v}^{n,M}_i(\z) \Bigg\} \Psi_{h}(\z) p(\z) d\z.
\end{split}
\ee
As before, the integrals are computed with Gaussian quadrature.
}
\paragraph{Approximation of quantities of interest.}
To analyze the effects of uncertainties on the kinetic model \eqref{eq:FP} we are generally interested in quantities of interest (QoI) $q[f]$ that can be computed from the distribution function. In the present setting, the QoI is simply the identity, i.e. $q[f] = f$. Anyway, more generally, it may be given by any moment of the kinetic solution.

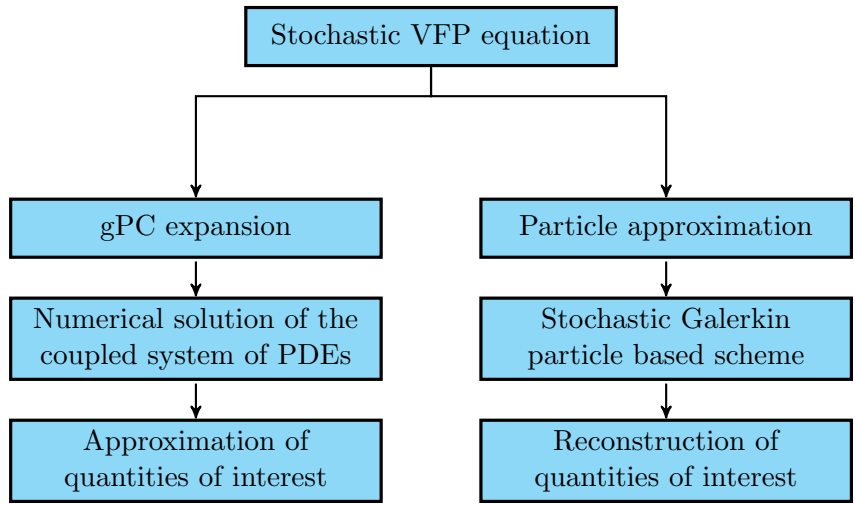
\begin{figure}[t]
\begin{center}
\begin{tikzpicture}
[node distance=0.5cm, start chain=going below,]
\node [punktchain] (Kinetic)  {\textnormal{Stochastic VFP equation}};
\node[below=2cm of Kinetic](dummy){};
\begin{scope}[start branch=venstre,every join/.style={->, thick, shorten <=1pt}, ]
\node [punktchain,left=of dummy] (SG-gPC)  {\textnormal{gPC expansion}};
\node [punktchain,join,on chain=going below] (sol SG)  {\textnormal{Numerical solution of the coupled system of PDEs}};
\node [punktchain,join, on chain=going below] (app)  {\textnormal{Approximation of quantities of interest}};
\end{scope}
\begin{scope}[start branch=hoejre,every join/.style={->, thick, shorten <=1pt}, ]
\node[punktchain,right=of dummy]  (SG-gPC-SDE)  {\textnormal{Particle approximation}};
\node [punktchain,join, on chain=going below] (MCgPC)  {\textnormal{Stochastic  Galerkin \\ particle based scheme  }};
\node [punktchain,join, on chain=going below] (app2)  {\textnormal{Reconstruction of quantities of interest }};
\end{scope}
\draw[|-,-|,->, thick,](Kinetic.south) |-+(0,-1em)-| (SG-gPC.north);
\draw[|-,-|,->, thick,] (Kinetic.south) |-+(0,-1em)-| (SG-gPC-SDE.north);
\end{tikzpicture}
\caption{\small{Numerical approaches to stochastic Vlasov-Fokker-Planck equations. The left branch describes the classical sG approach based on finite differences/volumes, whereas the right branch sketches the stochastic particle Galerkin approach.}}
\label{fig:scheme}
\end{center}
\end{figure}

In Figure \ref{fig:scheme} we sketch two approaches for the approximation of QoI of the Vlasov-Fokker-Planck equation. 
In the left branch we find the standard sG approach, in which we first consider the gPC approximation of the original problem that generates a coupled system of VFP-type PDEs. The resulting systems can be solved by classical deterministic finite difference/volume methods. In the right branch we describe instead the stochastic Galerkin particle reformulation of the problem, in which the gPC approximation is considered at the particle level.
We highlight that, in view of the above particle reformulation of the problem, the distribution is defined from a sample of $N$ particles $\{x_i^n(\z),v_i^n(\z)\}_{i=1}^N$ as follows
\[
{f_N^n(x,v,\z) = \dfrac{m(\z)}{N} \sum_{i=1}^N \delta(x-x_i^n(\z)) \otimes \delta(v-v_i^n(\z)),}
\]
{where $\delta(\cdot)$ is the Dirac delta}. In particular, we approximate this quantity as 
\[
{f_N^{n,M}(x,v,\z) = \dfrac{m(\z)}{N} \sum_{i=1}^N \delta(x-x_i^{n,M}(\z)) \otimes \delta(v-v_i^{n,M}(\z)),}
\]
having defined $\{x_i^{n,M}(\z),v_i^{n,M}(\z)\}_{i=1}^N$ in \eqref{eq:gPCE} and where $(\hat{x}^n_{i,h},\hat v^n_{i,h})$ are solution of \eqref{eq:sG_BGK} for the collision step and of \eqref{eq:sG_three_step1}-\eqref{eq:sG_three_step2}-\eqref{eq:sG_three_step3} for the transport step. 

\section{Numerical tests}
\label{sec:test}
In this section, we present several numerical results for the DSMC-sG Algorithm \ref{algorithm_sG_DSMC_BGK} and \ref{algorithm_sG_DSMC_transport}. In all the following tests, we consider the second order Strang splitting method exposed in Section \ref{sec:sG_particle_methods} on the time interval $[0,t_f]$. We restrict to bi-dimensional phase space $(x,v)\in\R\times\R$, for the {sake} of simplicity, however the results hold also in higher-dimensional phase space. We reconstruct the densities through {histograms} with $N_\ell=100$ cells in the spatial domain and $N_v=200$ cells in the velocity domain. The Poisson equation is solved with a finite difference method on {a uniform grid} of $N_x=N_\ell+1=101$ points, with periodic or Dirichelet boundary conditions  {\cite{Sonnendrucker2013}}, adopting a stochastic collocation approach^^>\cite{ xiu2010}. In all the test we considered uniform uncertainties, that are linked to Legendre polynomials in the Wiener-Askey scheme, see^^>\cite{xiu2002}. {The details of the initial sampling techniques are presented in Appendix \ref{appendix}.}

First we consider the linear and nonlinear Landau damping, then we check the spectral convergence of the scheme and we investigate another standard plasma test, the two stream instability, both in linear and nonlinear case. In the end, we investigate the hydrodynamic limit of the model with the standard Sod shock tube test.


\subsection{Test 1: Landau damping}
The Landau damping^^>\cite{landau1965,Villani2013} in one of the greatest theoretical {results} of the collisionless plasma physics. It corresponds to the exponential decay of the electromagnetic waves without energy dissipation, as a result of the wave-particle interaction. The charged particles in motion with a velocity lower than the phase velocity of the wave gain energy from it, leading to the electromagnetic energy damping. On the contrary, particles with a velocity greater than the phase velocity yield energy to the wave, causing its energy increasing. As a result, although the total energy of the Vlasov-Poisson equation is conserved for every $\z\in\R^{d_\z}$, i.e.
\[
\dfrac{d}{dt} \left[ \dfrac{1}{2}\int_{\R^3}\int_{\R^3}  |v|^2 f(x,v,t,\z) dx dv + \dfrac{1}{2}\int_{\R^3} |E(x,t,\z)|^2 dx \right] = 0,
\]    
if there are more particles with $v<v_\phi$ than particles with $v>v_\phi$, we may observe the decay of the L$^2$-norm of the electric field 
\be \label{eq:el_energy}
\mathcal{E}(t,\z) = \left( \int_{\R^3} |E(x,t,\z)|^2 dx \right)^{\frac{1}{2}}
\ee
with a specific damping rate $\gamma$. If we consider also a small but finite collision frequency $\nu$, the damping is balanced by the interactions among the particles and, in the limit $\nu\to+\infty$, we recover the fluid limit of the Vlasov equation with the electric field fluctuating around a constant value. 

To observe the Landau damping, we consider a wave perturbation of the local Maxwellian distribution. If the perturbation is small, we are in the so-called \textit{linear Landau damping} regime, if the wave amplitude increases, we get the \textit{nonlinear Landau damping} regime. 

\paragraph{Linear case.}
We consider, at the initial time, an uncertain perturbation of the local equilibrium 
\be \label{eq:init_linear}
f_0(x,v,\z) = \left( 1 + \alpha(\z) \cos(k x) \right) \dfrac{1}{\sqrt{2\pi T}} e^{-\frac{v^2}{2T}},
\ee
with $x\in[0,2\pi/k]$, $k$ wave number and $\alpha(\z)$ small random amplitude. The temperature is fixed $T=1$. 

To observe the phenomenon, the velocity domain has to be larger than the phase velocity of the wave $v_\phi = \omega / k$, where $\omega$ is the frequency and $v_\phi$ is the singularity of the dispersion relation. According to^^>\cite{Crouseilles2004}, the frequency $\omega$ can be estimated as
\[
\omega^2 = 1 + 3k^2,
\]
while, for small wave numbers, a more accurate formula^^>\cite{McKinstrie1999} gives
\[
\omega^2 = 1 + 3k^2 + 6 k^4 + 12 k^6.
\]

In the collisionless regime^^>\cite{Chen1974}, we are able to provide an approximation also for the damping rate $\gamma$, that is
\[
\gamma = -\sqrt{\dfrac{\pi}{8}} \dfrac{1}{k^3} \exp\left( -\dfrac{1}{2 k^2} - \dfrac{3}{2}\right)
\]
for large $k$, and  
\[
\gamma = -\sqrt{\dfrac{\pi}{8}} \left( \dfrac{1}{k^3} - 6 k \right) \exp\left( -\dfrac{1}{2 k^2} - \dfrac{3}{2} - 3k^2 - 12k^4 \right)
\]
for smaller wave number. 

We consider $N=10^7$ particles, $\Delta t=0.1$, $M=5$ for the stochastic Galerkin projection and $\alpha(\z) = \frac{1}{20}+\frac{1}{10}\z$, $\z\sim\mathcal{U}([0,1])$ in \eqref{eq:init_linear}. We select $k=0.5$ and we investigate different collisional {regimes} corresponding to the choices $\nu=0,1,10^3$. We apply periodic boundary conditions on the particles in the space domain according to \eqref{eq:pbc} and on the Poisson equation.

\begin{figure}
	\centering
	\includegraphics[width = 0.325\linewidth]{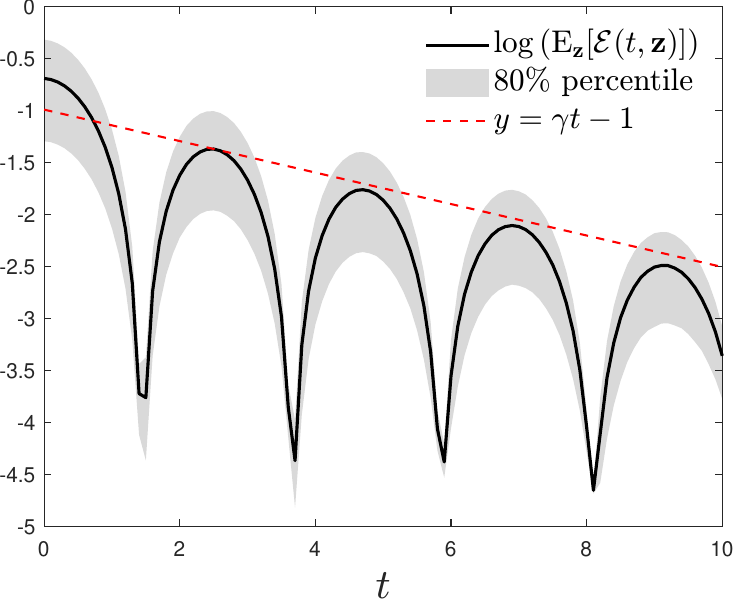}
	\includegraphics[width = 0.325\linewidth]{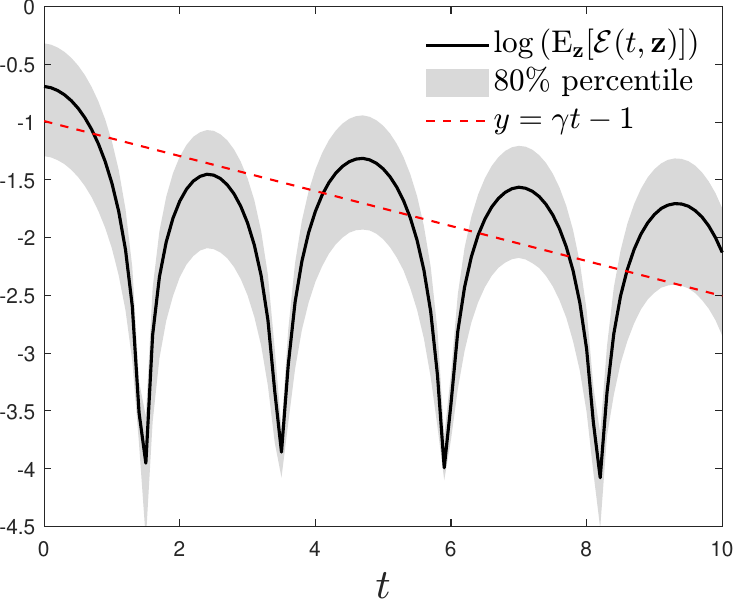}
	\includegraphics[width = 0.325\linewidth]{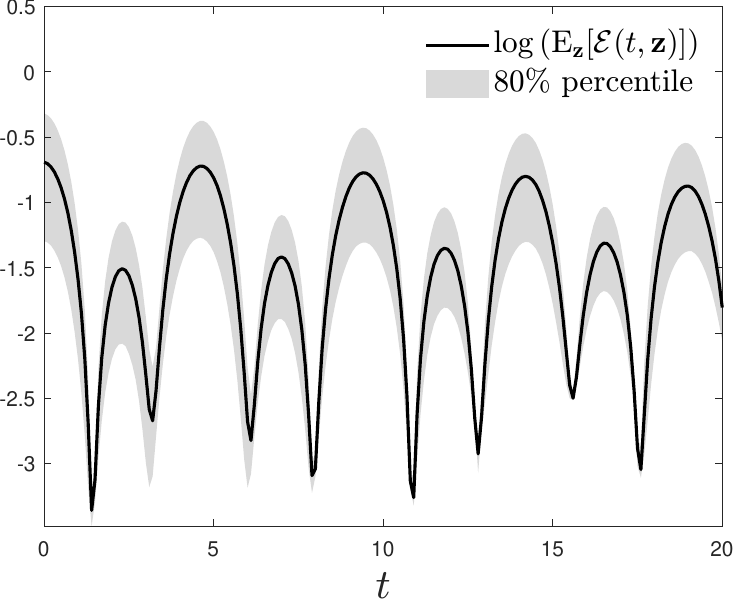}
	\includegraphics[width = 0.325\linewidth]{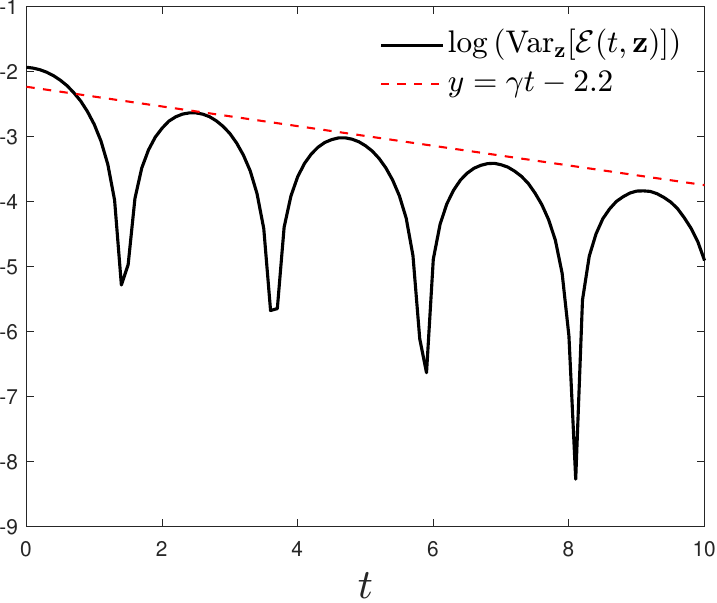} 
	\includegraphics[width = 0.325\linewidth]{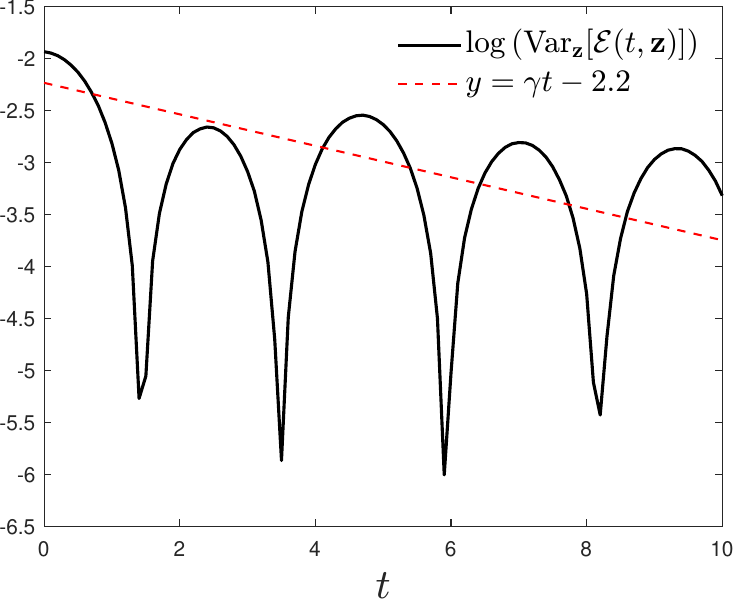} 
	\includegraphics[width = 0.325\linewidth]{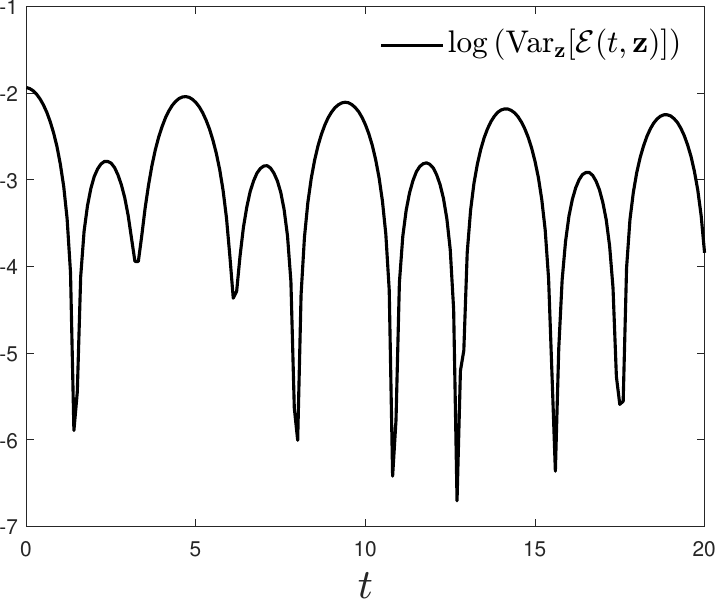}
	\caption{\small{\textbf{Test 1: linear Landau damping}. Logarithm of $\mathbb{E}_{\z}[\mathcal{E}(t,\z)]$ (top row) and $\textrm{Var}_{\z}[\mathcal{E}(t,\z)]$ (bottom row) for $\nu=0$ (left), $\nu=1$ (centre) and $\nu=10^3$ (right) in linear Landau damping. We choose $k=0.5$, $N= 10^7$ particles, $M=5$ and $\Delta t = 0.1$. We considered the initial condition \eqref{eq:init_linear} with $\alpha(\z) = \frac{1}{20}+\frac{1}{10}\z$, $\z\sim\mathcal{U}([0,1])$. The theoretical damping rate is $\gamma=-0.1533$.}}
	\label{fig:test_2_linear}
\end{figure}

In Figure \ref{fig:test_2_linear} we show the logarithm of the expectation and the variance of $\mathcal{E}(t,\z)$. As expected, in the collisionless scenario $\nu=0$, the numerical results are in agreement with the theoretical damping rate $\gamma=-0.1533$ and with the approximation of the electric field obtained, e.g., in^^>\cite{Sonnendrucker2013}. With a collision frequency $\nu=1$, we observe that the interactions between the particle balance the damping and the electric field exhibits a slower decreasing with respect to the collisionless case. With a bigger frequency $\nu=10^3$, we note that the electric field fluctuates around a constant value, as expected (see, e.g.,^^>\cite{Dimarco2015}). 

As we can observe from the bottom row of Figure \ref{fig:test_2_linear}, the dynamics of the variance shows a time evolution that is similar to the expected value. This is an effect of the linear perturbation regime, since the oscillation frequency of the electromagnetic wave depends only on the wave number.

\paragraph{Nonlinear case.}
We investigate now the nonlinear case, that is, we consider the initial distribution \eqref{eq:init_linear} with  $N= 5\cdot 10^7$ particles and the same computational setting as before with the only difference of a greater perturbation amplitude. In particular, we choose $\alpha(\z) = \frac{2}{5}+\frac{3}{5}\z$, $\z\sim\mathcal{U}([0,1])$, in the collisionless scenario $\nu=0$, and $\alpha(\z) = \frac{1}{5}+\frac{2}{5}\z$, $\z\sim\mathcal{U}([0,1])$, for $\nu=1, 10^3$.

The analytical estimations obtained in the previous paragraph via the linear theory no longer hold. Since the electric potential is larger, it could happen that particles remain trapped in the potential well of the electromagnetic wave. This is the so-called \textit{electron trapping}. The trapped particles, fluctuating into the well, cause both damping and growth of the electric field. In literature there exist several approximations for the damping $\gamma_d$ and growth $\gamma_g$ rates in the nonlinear regime, from which we can compare our numerical results, see e.g.^^>\cite{Campospinto2014, chacon2016, filbet2003, rossmanith2011, liu2017, Dimarco2015} and the references therein.

\begin{figure}
	\centering
	\includegraphics[width = 0.325\linewidth]{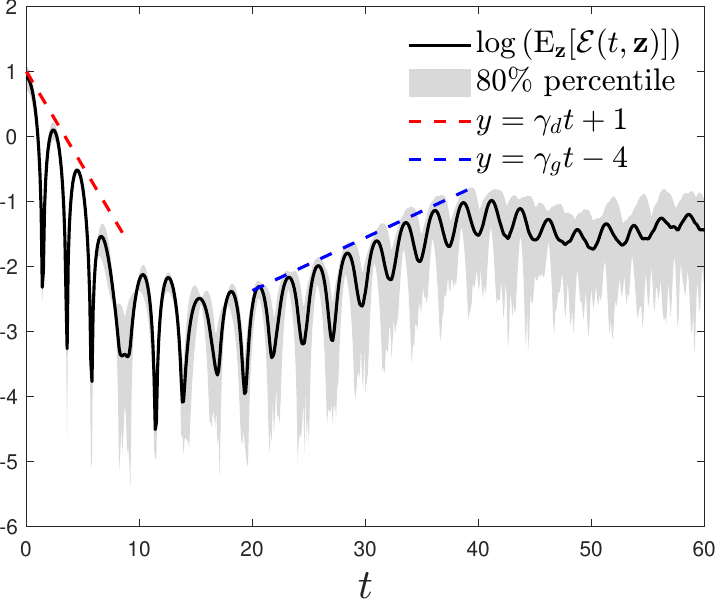}
	\includegraphics[width = 0.325\linewidth]{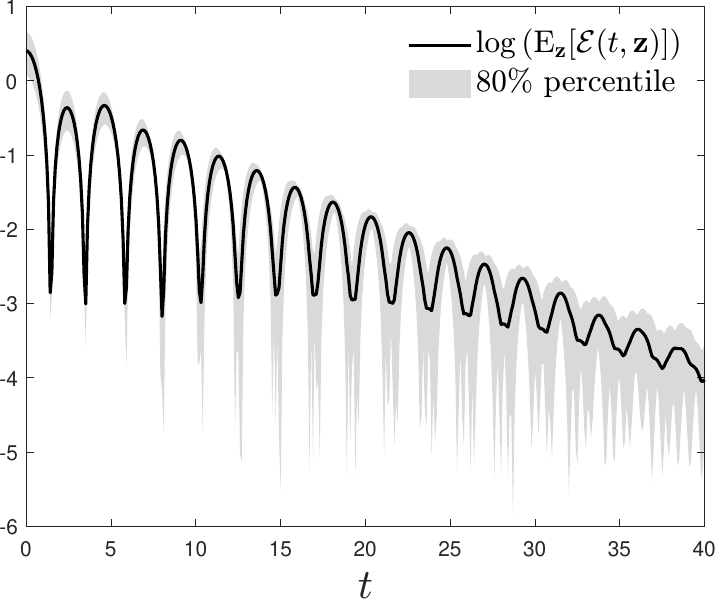}
	\includegraphics[width = 0.325\linewidth]{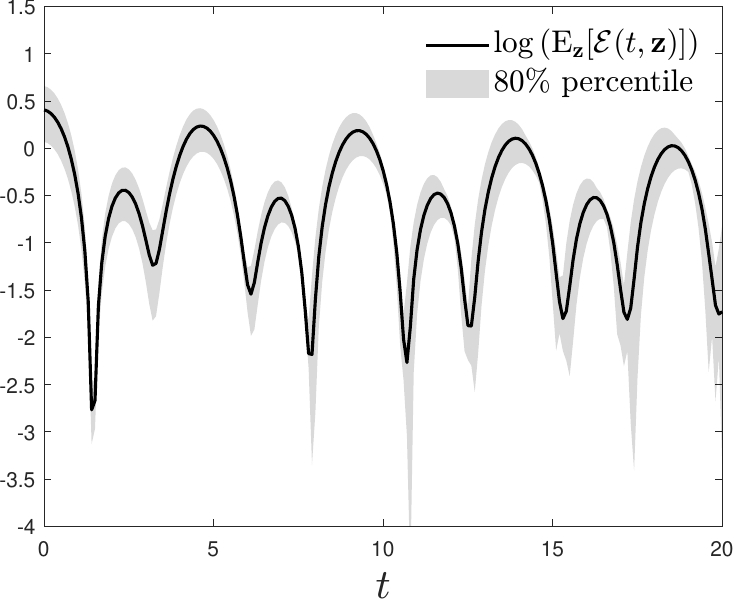}
	\includegraphics[width = 0.325\linewidth]{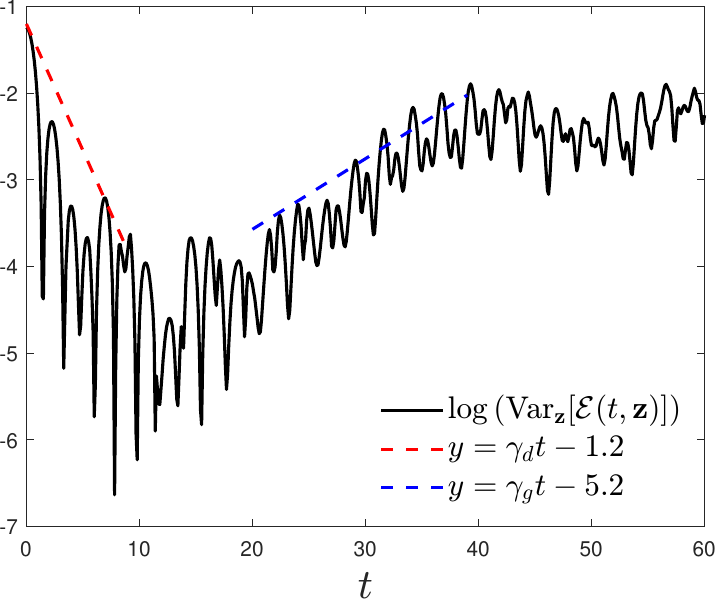} 
	\includegraphics[width = 0.325\linewidth]{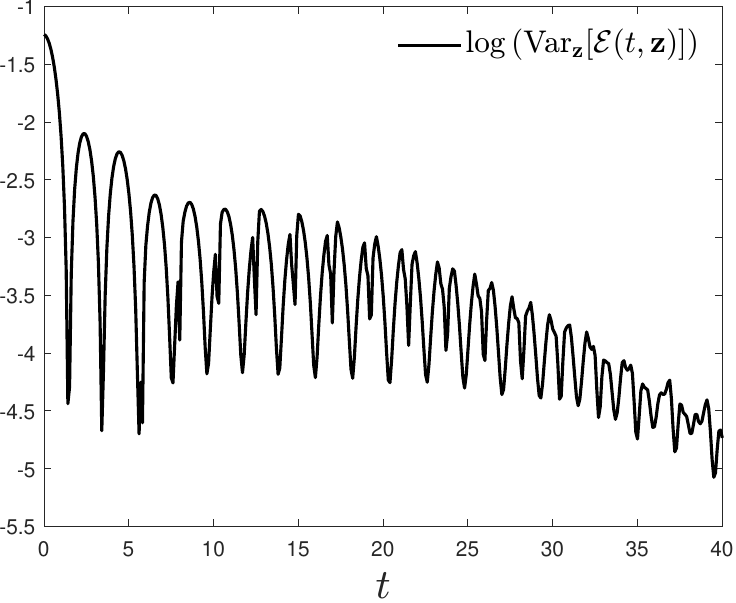} 
	\includegraphics[width = 0.325\linewidth]{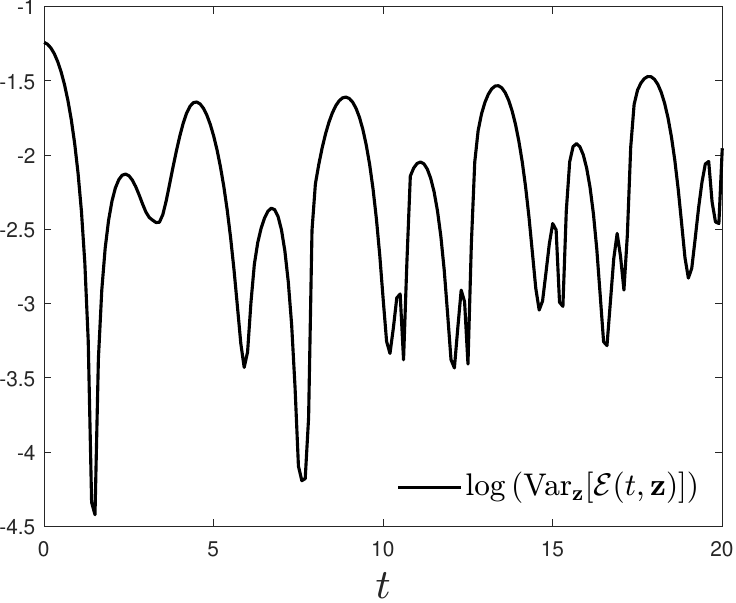}
	\caption{\small{\textbf{Test 1: nonlinear Landau damping}. Logarithm of $\mathbb{E}_{\z}[\mathcal{E}(t,\z)]$ (top row) and $\textrm{Var}_{\z}[\mathcal{E}(t,\z)]$ (bottom row) for $\nu=0$ (left), $\nu=1$ (centre) and $\nu=10^3$ (right) in nonlinear Landau damping. We choose $k=0.5$, $N= 5\cdot 10^7$ particles, $M=5$ and $\Delta t = 0.1$. We considered the initial condition \eqref{eq:init_linear} with $\alpha(\z) = \frac{2}{5}+\frac{3}{5}\z$, $\z\sim\mathcal{U}([0,1])$. In the collisionless scenario the damping rate is $\gamma_d=-0.2920$ and the growth rate is $\gamma_g=0.0815$.}}
	\label{fig:test_2_non_linear}
\end{figure}

As we can observe in Figure \ref{fig:test_2_non_linear}, in the collisionless scenario the numerical results are in accordance with the theoretical approximations of the rates, i.e. $\gamma_d=-0.2920$ and $\gamma_g=0.0815$. With collision frequencies $\nu=1,10^3$, we observe that the interactions between the particles induce slight different behaviours of the electric energy, balancing the damping and growth behaviour. We observe also that the variance, as the time increases, exhibits a shift in the peaks of oscillation. Unlike the linear case, the frequency of the electric wave appears as dependent on the initial uncertain amplitude.

\begin{figure}
	\centering
	\includegraphics[width = 0.325\linewidth]{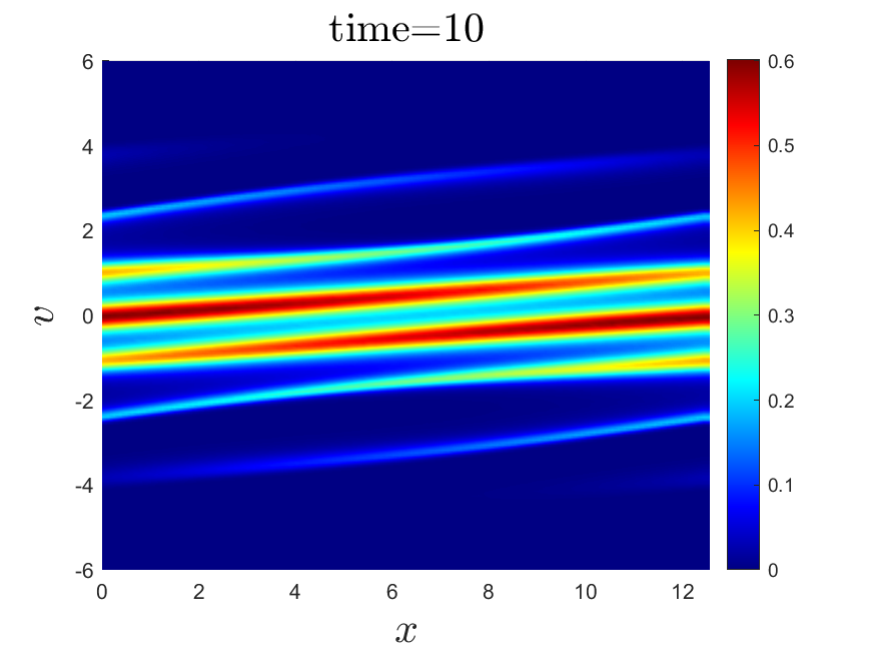}
	\includegraphics[width = 0.325\linewidth]{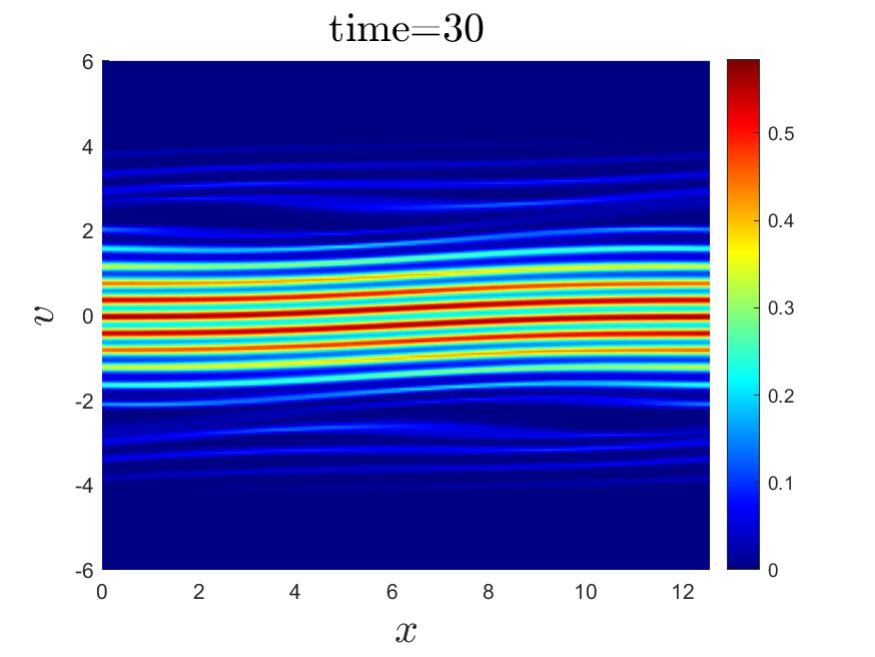}
	\includegraphics[width = 0.325\linewidth]{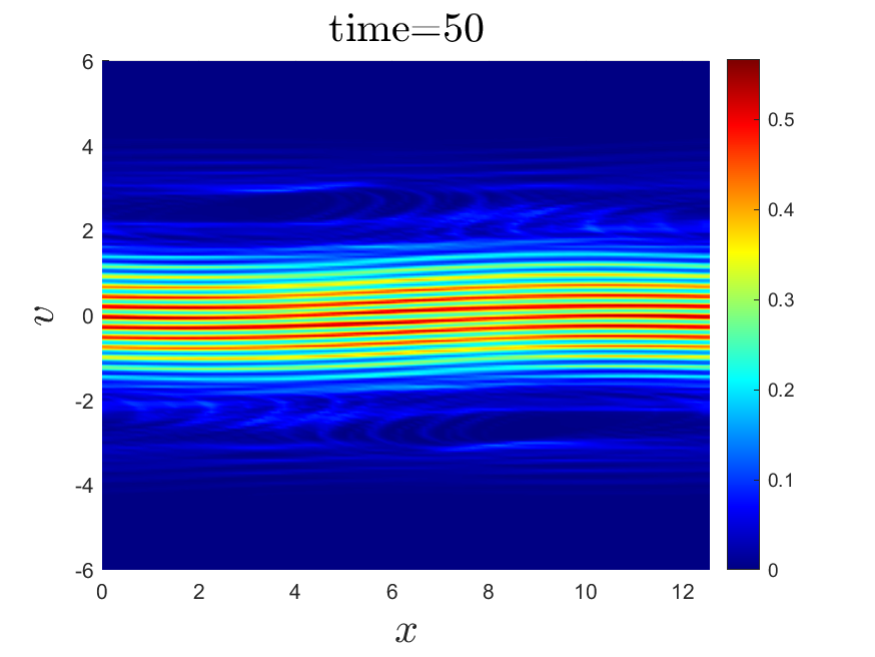}
	\includegraphics[width = 0.325\linewidth]{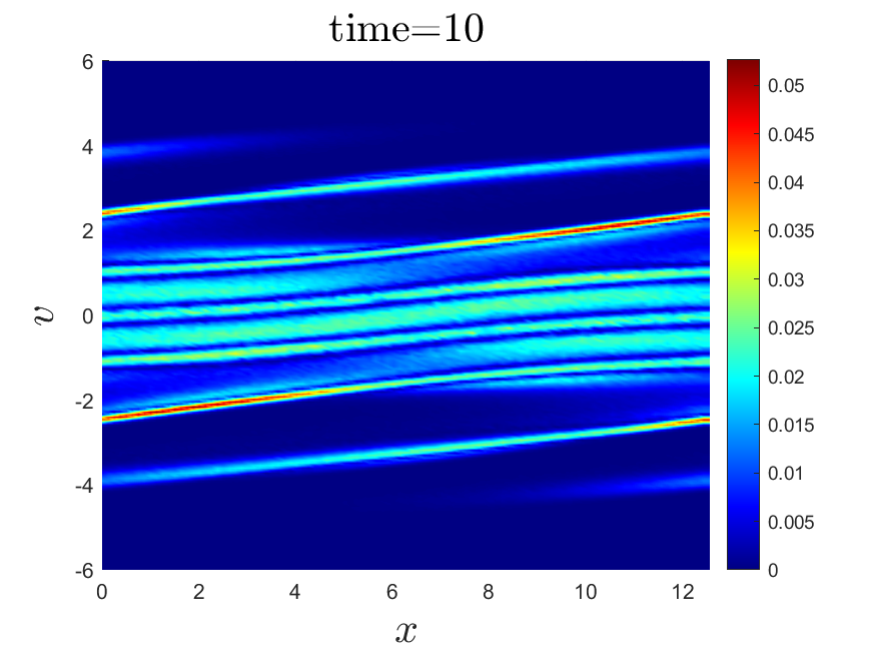} 
	\includegraphics[width = 0.325\linewidth]{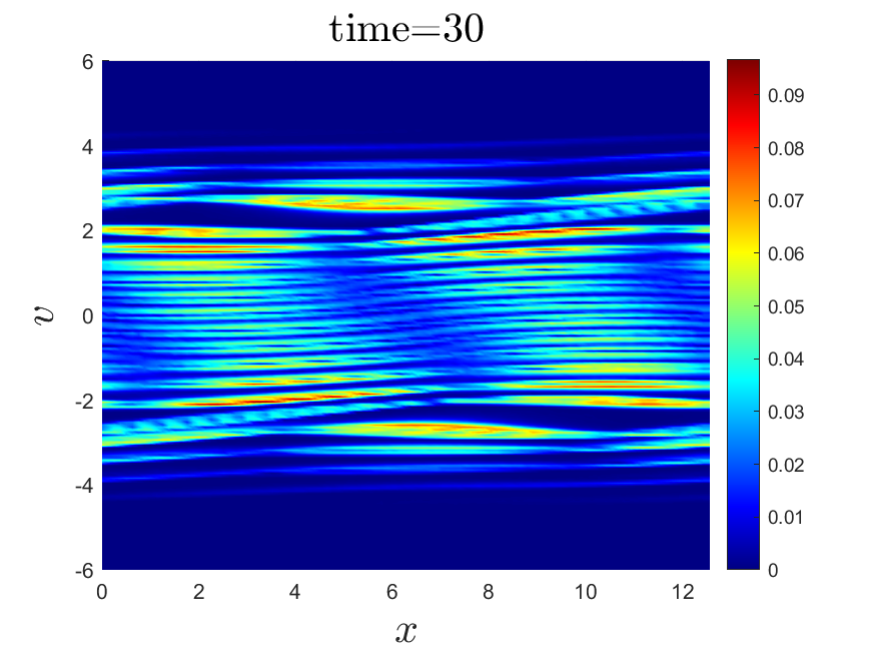} 
	\includegraphics[width = 0.325\linewidth]{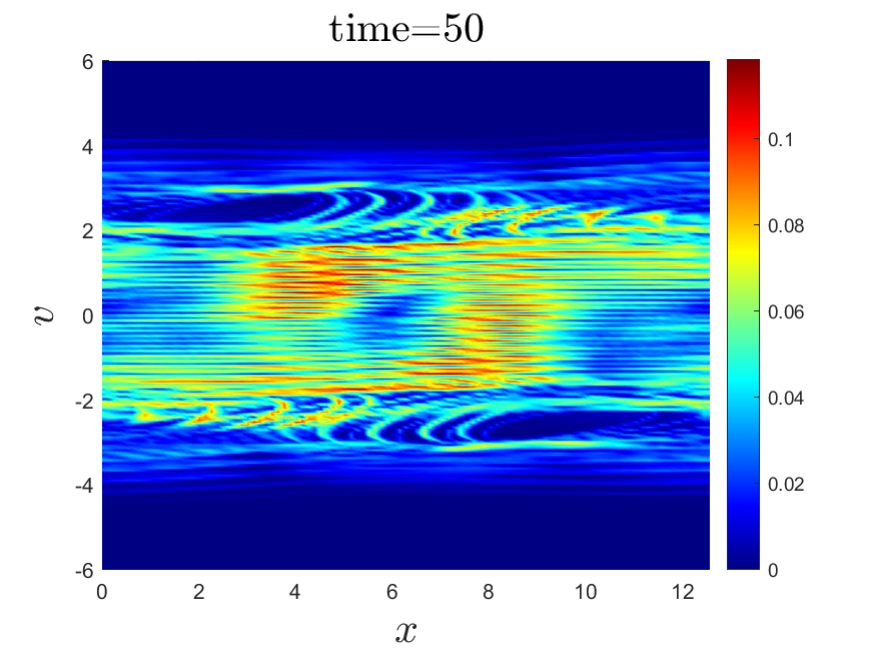}
	\caption{\small{\textbf{Test 1: nonlinear Landau damping}. We report $\mathbb{E}_{\z}[f(x,v,t,\z)]$ (top row) and $\textrm{Var}_{\z}[f(x,v,t,\z)]$ (bottom row) at fixed times $t=10,30,50$, in collisionless nonlinear Landau damping. We choose $k=0.5$, $N=5 \cdot 10^7$ particles, $M=5$ and $\Delta t=0.1$. We considered the initial condition \eqref{eq:init_linear} with $\alpha(\z) = \frac{2}{5}+\frac{3}{5}\z$, $\z\sim\mathcal{U}([0,1])$.}}
	\label{fig:test_2_non_linear_distributions}
\end{figure}

In Figure \ref{fig:test_2_non_linear_distributions} we show the typical profiles of the expectation and variance of the distributions in the collisionless case, at fixed times $t=10,30,50$ and we note that the expected value and the variance in $\z$ exhibit similar foliation patterns.

\subsection{Test 2: spectral convergence}
We numerically check for the convergence in the space of the random parameters of the Algorithm \ref{algorithm_sG_DSMC_transport}, i.e. the collisionless Vlasov-Poisson system with $\nu=0$. 

We consider the initial distribution
\be \label{eq:init_conv}
f_0(x,v,\z) = \rho(x) \dfrac{1}{\sqrt{2\pi T(\z)}} e^{-\frac{v^2}{2T(\z)}}, 
\ee
with initial random temperature $T(\z) = \frac{4}{5} + \frac{2}{5} \z$, $\z\sim\mathcal{U}([0,1])$ and mass
\[
\rho(x) = \dfrac{1}{\sqrt{\pi}} e^{-(x-6)^2}
\]
in the physical domain $x\in[0,4\pi]$. We use $N=10^6$ particles and $\Delta t=0.1$, solving the Poisson equation with periodic boundary conditions.

\begin{figure}
	\centering
	\includegraphics[width = 0.5\linewidth]{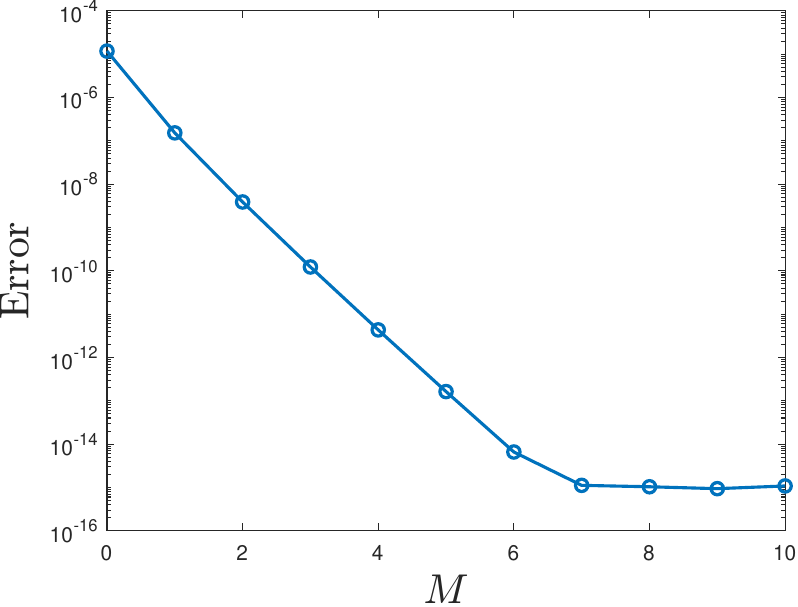}
	\caption{\small{\textbf{Test 2: spectral convergence}. $L^2$ error in the evaluation of $\mathcal{E}(t,\z)$ of the sG particle scheme at fixed time $t=1$. We consider $N=10^6$ particles, $\Delta t=0.1$ and a reference solution computed with $M=30$. We consider the initial condition \eqref{eq:init_conv} with $T(\z) = \frac{4}{5} + \frac{2}{5} \z$, $\z\sim\mathcal{U}([0,1])$.}}
	\label{fig:test_1_spectral}
\end{figure}

To observe the convergence, we consider a reference sG particle solution obtained with $M=30$ and we store the initial data, the only source of randomness of the scheme, so that the same data is used for the different values of $M$.
In Figure \ref{fig:test_1_spectral} we show in semilogarithmic scale the decay of the $L^2$ error in the evaluation of the observable $\mathcal{E}(t,\z)$ defined in \eqref{eq:el_energy}, for increasing $M$ and at fixed time $t=1$. Exponential decay of the error up to machine precision is observed, which shows that the sG particle scheme is spectrally accurate in the space of the random parameters.

\subsection{Test 3: two stream instability}
Another classical effect of the collisionless plasma is the so-called two stream instability. If we initialize the electron velocities as a non-isotropic two beam distribution, contrary to the Landau damping, we have more particles with $v>v_\phi$ than particles with $v<v_\phi$. As a consequence, the electromagnetic wave gains energy from the electrons, giving birth to instabilities represented by typical swirling and warped pattern of the particle distribution function. In the following, we will investigate the linear and the nonlinear case.

\paragraph{Linear case.}
We consider the initial distribution
\be \label{eq:init_linear_twostream}
f_0(x,v,\z) = \left( 1 + \alpha(\z) \cos(k x) \right) \dfrac{1}{\sqrt{2\pi T}} \left( e^{-\frac{(v-\bar{v})^2}{2T}} + e^{-\frac{(v+\bar{v})^2}{2T}}\right),
\ee
with fixed temperature $T=1$ and $\bar{v}=2.4$. As before, we have $x\in[0,2\pi/k]$.

In the collisionless scenario, if the perturbation amplitude is small enough, after a certain amount of time the logarithm of the L$^2$-norm of the electric energy grows linearly with a specific rate $\gamma$, see e.g.^^>\cite{Chen1974, Xiao2021, liu2017}. 

\begin{figure}
\centering
\includegraphics[width = 0.45\linewidth]{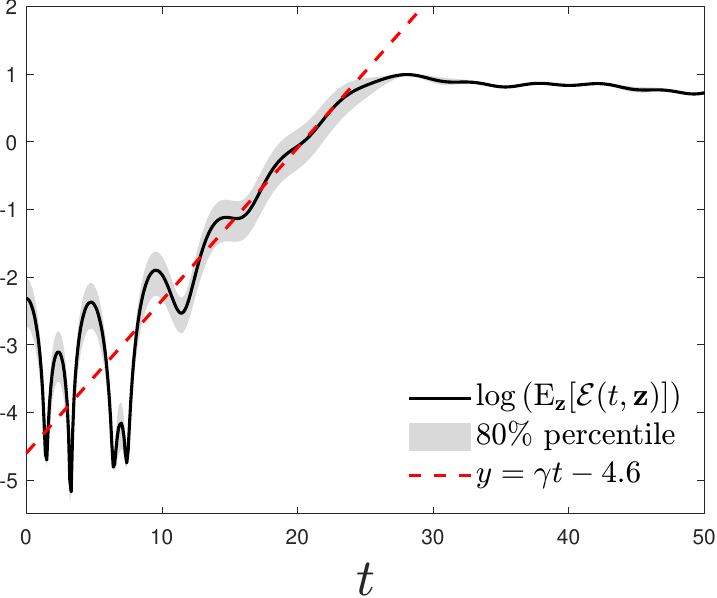}
\includegraphics[width = 0.45\linewidth]{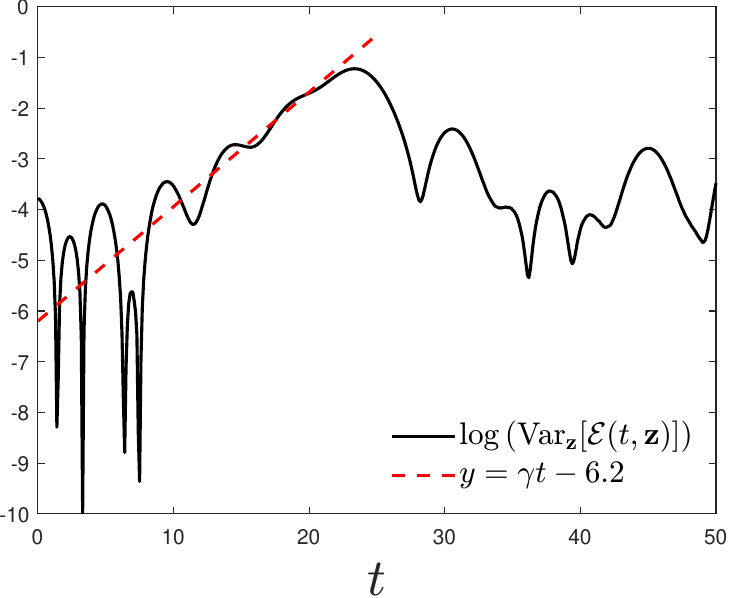}
\caption{\small{\textbf{Test 3: linear two stream instability}. Logarithm of $\mathbb{E}_{\z}[\mathcal{E}(t,\z)]$ (left) and $\textrm{Var}_{\z}[\mathcal{E}(t,\z)]$ (right) in collisionless linear two stream instability. We choose $k=0.2$, $N= 5 \cdot 10^7$ particles, $M=5$ and $\Delta t = 0.1$. We considered the initial condition \eqref{eq:init_linear_twostream} with $\alpha(\z) = 3\cdot 10^{-3} + 4\cdot 10^{-3}\z$, $\z\sim\mathcal{U}([0,1])$. The theoretical growth rate is $\gamma=0.2258$.}}

\label{fig:test_3_linear_twostreams}
\end{figure}

We choose $N=5\cdot 10^7$ particles $M=5$, $k=0.2$, $\nu=0$ and $\alpha(\z) = 3\cdot 10^{-3} + 4\cdot 10^{-3}\z$, $\z\sim\mathcal{U}([0,1])$. As in the Landau damping tests, we apply periodic boundary conditions on the particles and on the Poisson equation. In Figure \ref{fig:test_3_linear_twostreams}, we show the expectation and variance of the electric energy and we observe that the results fit the theoretical linear growth of $\gamma=0.2258$. After a certain point, we note that the mean reaches a stable value with a non zero variance, meaning that the uncertainties are not dampened by the dynamics of the model. 

\begin{figure}
\centering
\includegraphics[width = 0.325\linewidth]{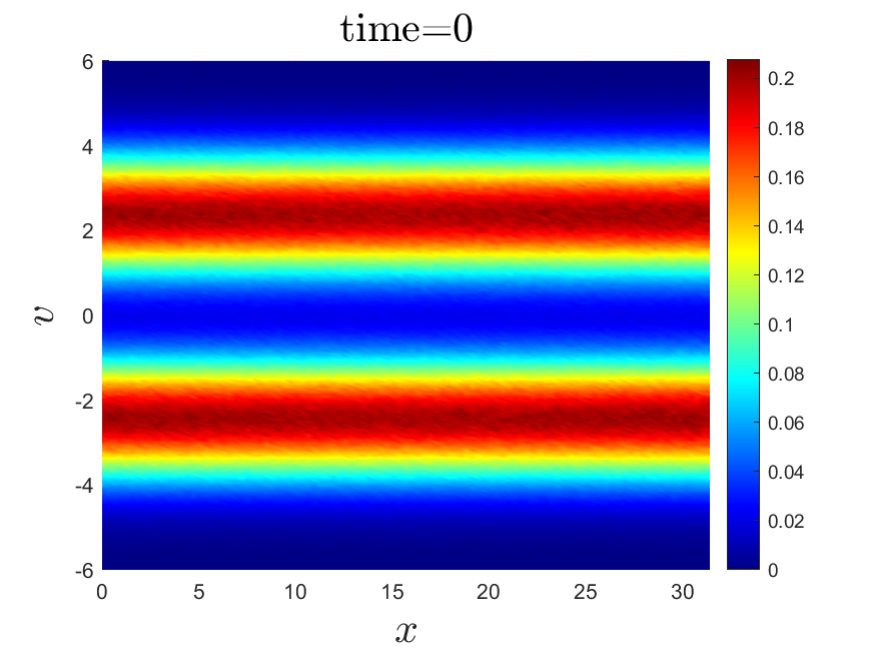}
\includegraphics[width = 0.325\linewidth]{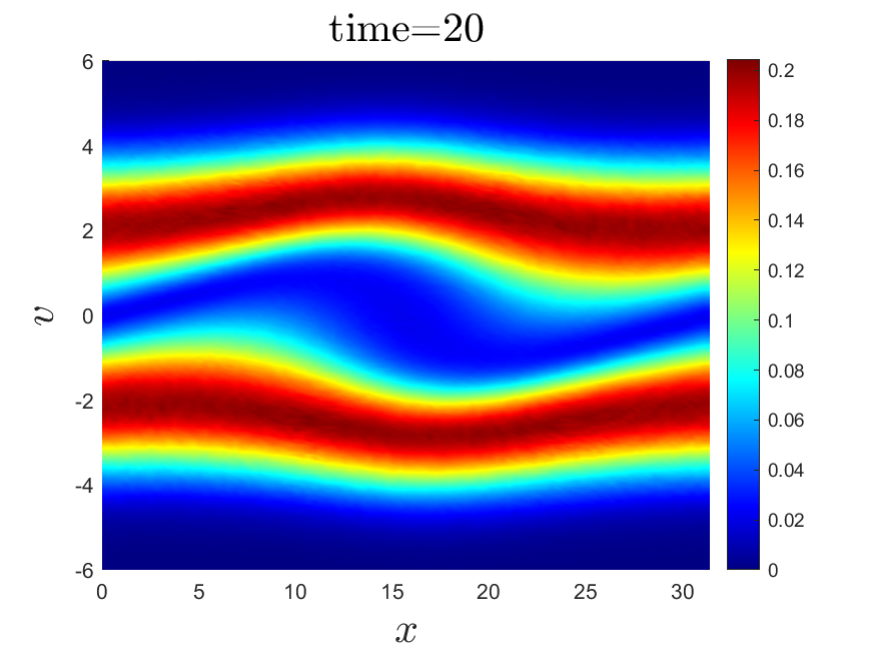}
\includegraphics[width = 0.325\linewidth]{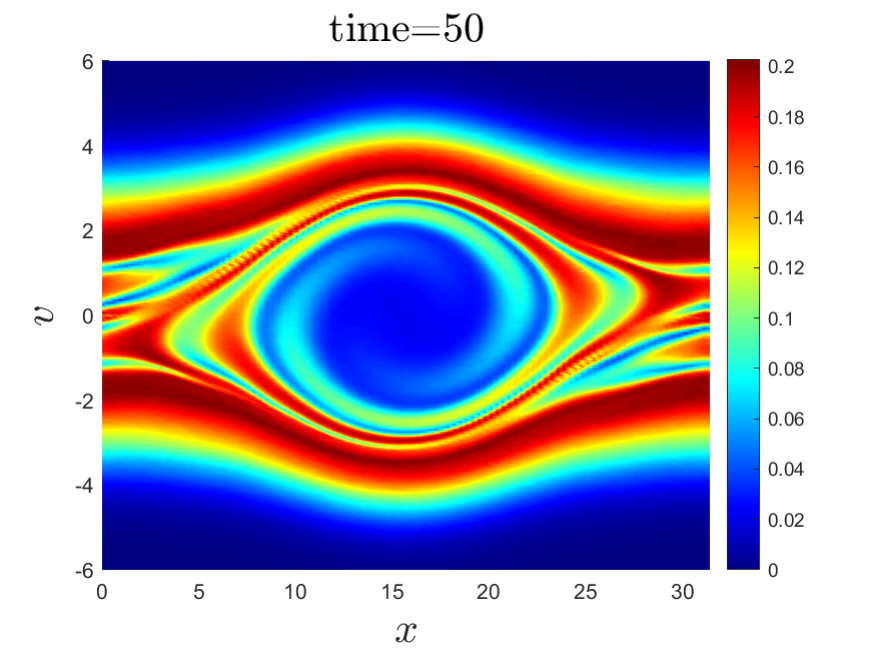}
\includegraphics[width = 0.325\linewidth]{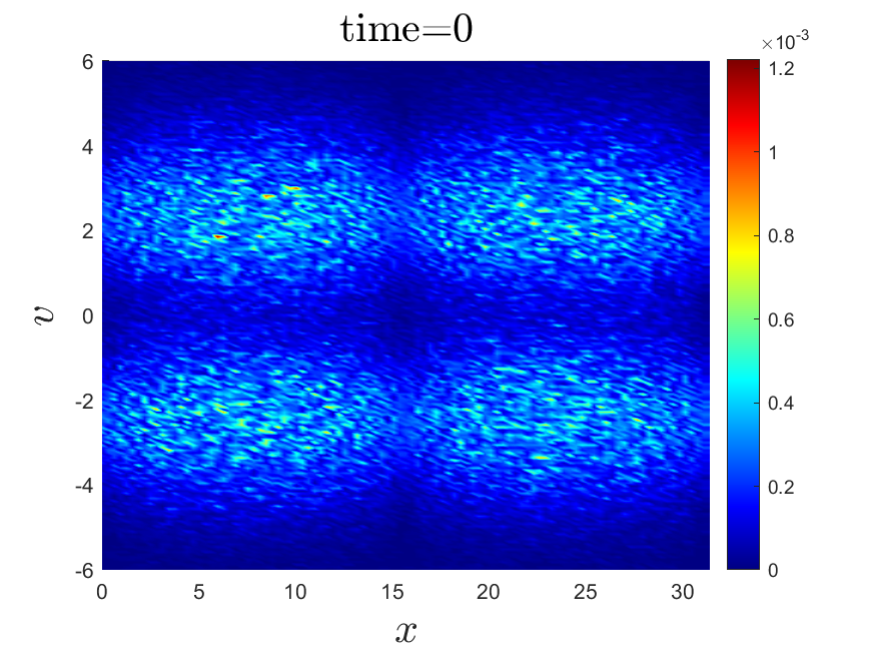} 
\includegraphics[width = 0.325\linewidth]{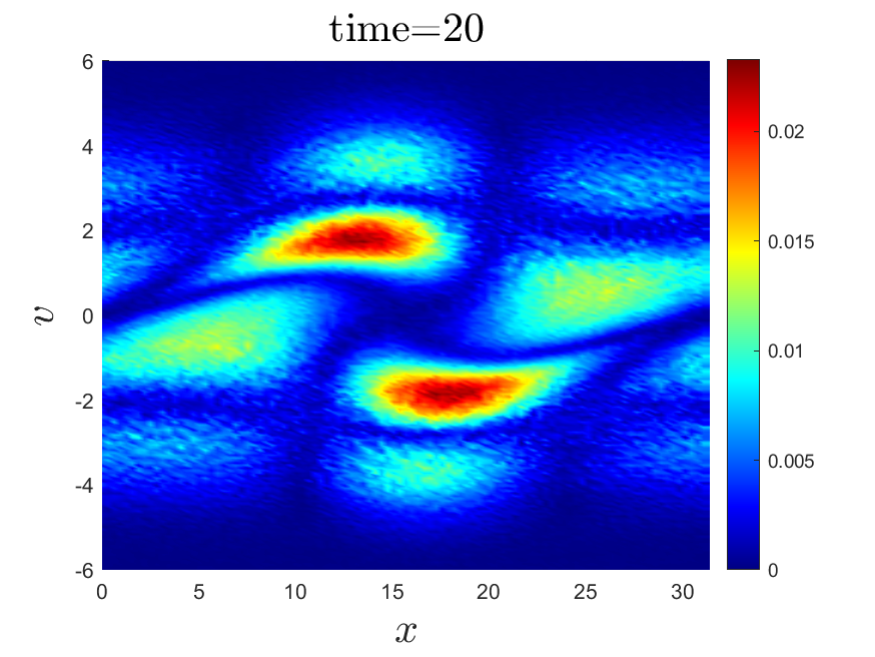} 
\includegraphics[width = 0.325\linewidth]{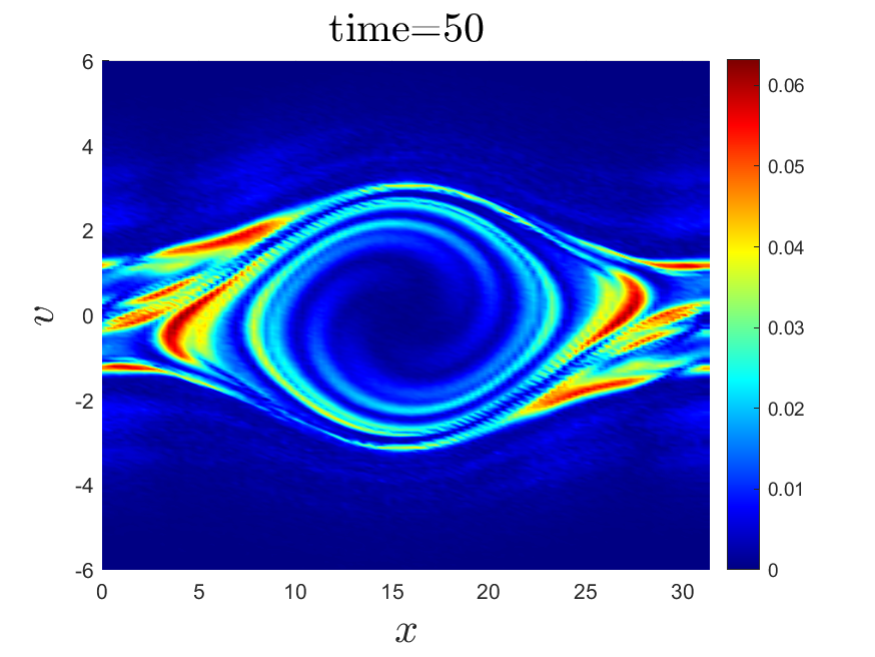}
\caption{\small{\textbf{Test 3: linear two stream instability}. We report $\mathbb{E}_{\z}[f(x,v,t,\z)]$  (top row) and $\textrm{Var}_{\z}[f(x,v,t,\z)]$ (bottom row) at fixed times $t=0,20,50$, in collisionless linear two stream instability. We choose $k=0.2$, $N=5 \cdot 10^7$ particles, $M=5$ and $\Delta t=0.1$. We considered the initial condition \eqref{eq:init_linear_twostream} with $\alpha(\z) = 3\cdot 10^{-3} + 4\cdot 10^{-3}\z$, $\z\sim\mathcal{U}([0,1])$.}}
\label{fig:test_3_linear_twostreams_distributions}
\end{figure}

In Figure \ref{fig:test_3_linear_twostreams_distributions}, we present the expectations and variances of the distribution at fixed times $t=0,20,50$. As the time increases, the typical whirling profiles of the linear two stream instability become more evident both in the mean and the variance. 

\paragraph{Nonlinear case.}
We consider here the nonlinear two stream instability. We choose the initial conditions 
\be \label{eq:init_nonlinear_twostream}
f_0(x,v,\z) = \left( 1 + \alpha(\z) \cos(k x) \right) \dfrac{1}{\sqrt{2\pi T}} \left( e^{-\frac{(v-\bar{v})^2}{2T}} + e^{-\frac{(v+\bar{v})^2}{2T}}\right),
\ee
with fixed temperature $T=0.3$, $\bar{v}=0.99$ and wave number $k=2/13$. The physical domain is $x\in[0,2\pi/k]$. We apply periodic boundary conditions on the particles and on the Poisson equation.  

We choose $N=5\cdot 10^7$ particles, $M=5$, $\Delta t=0.1$ and the perturbation amplitude $\alpha(\z) = 4\cdot 10^{-2} + 2\cdot 10^{-2}\z$, $\z\sim\mathcal{U}([0,1])$. In the following, we investigate different collisional \textcolor{magenta}{regimes} corresponding to the choices $\nu=0,1,10^3$. 

\begin{figure}
	\centering
	\includegraphics[width = 0.325\linewidth]{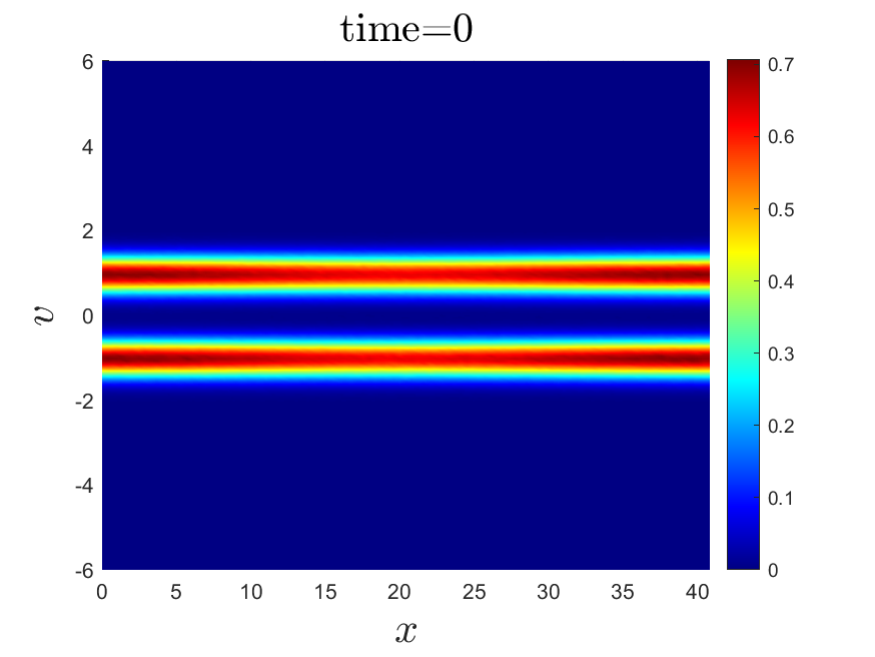}
	\includegraphics[width = 0.325\linewidth]{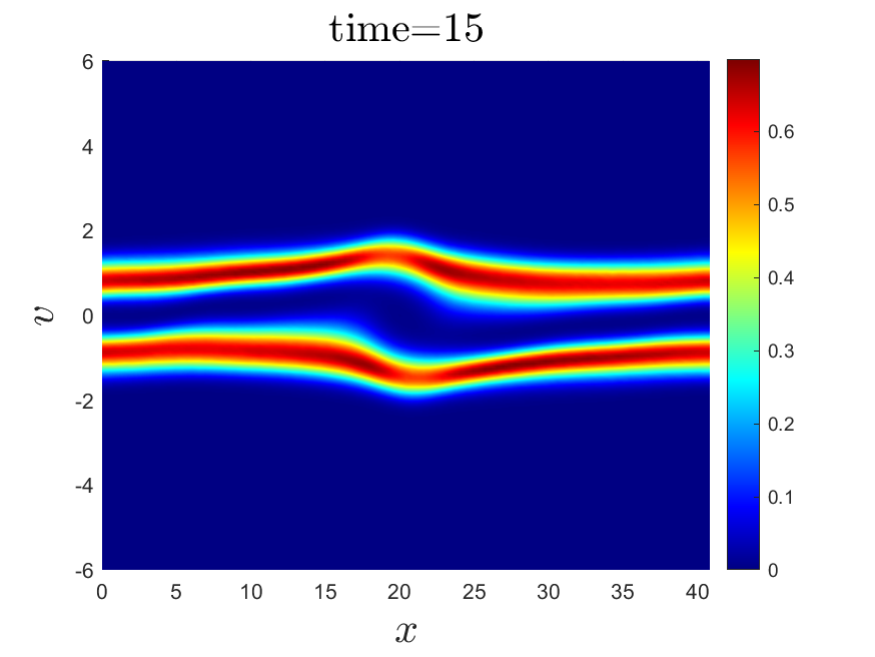}
	\includegraphics[width = 0.325\linewidth]{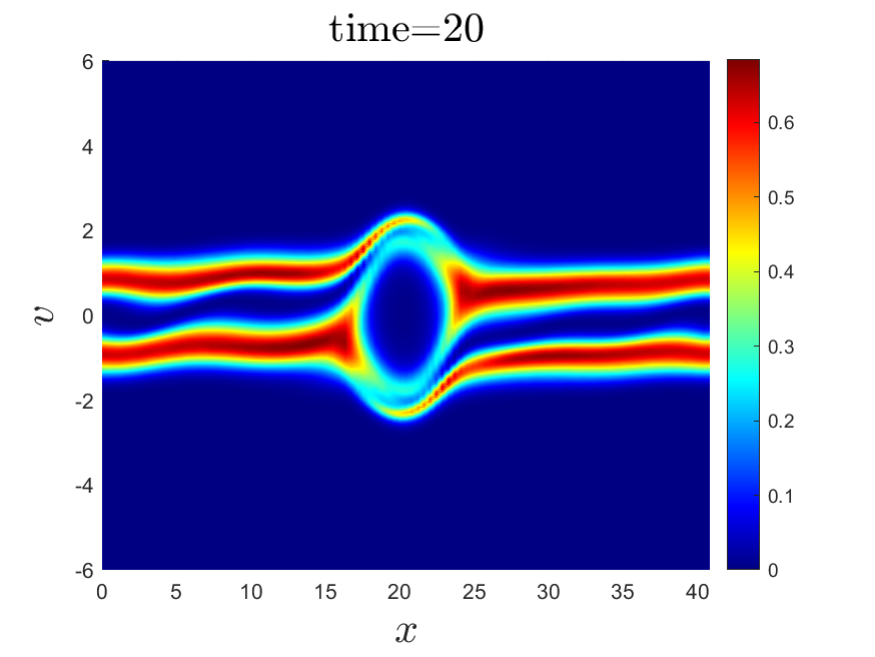}
	\includegraphics[width = 0.325\linewidth]{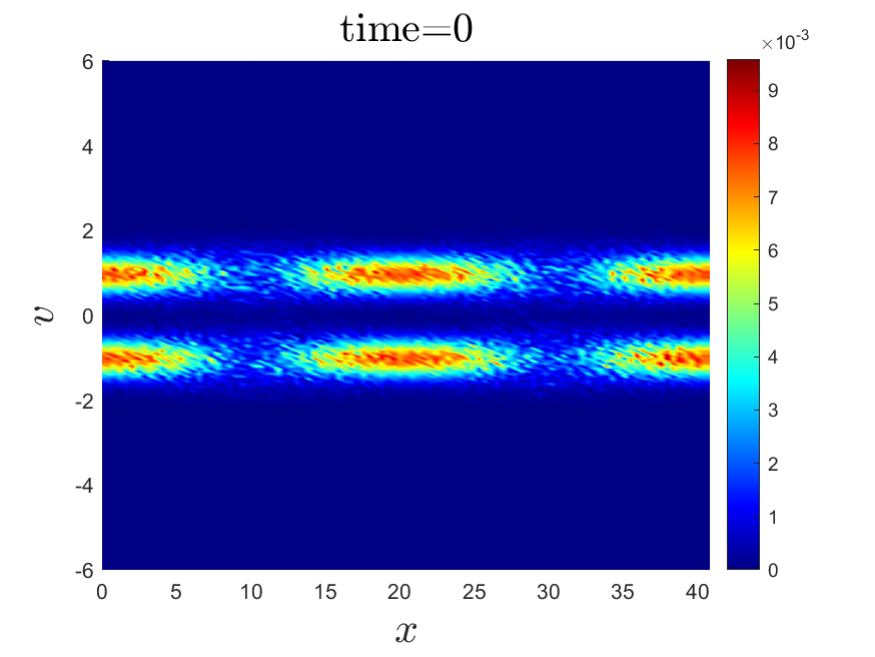}
	\includegraphics[width = 0.325\linewidth]{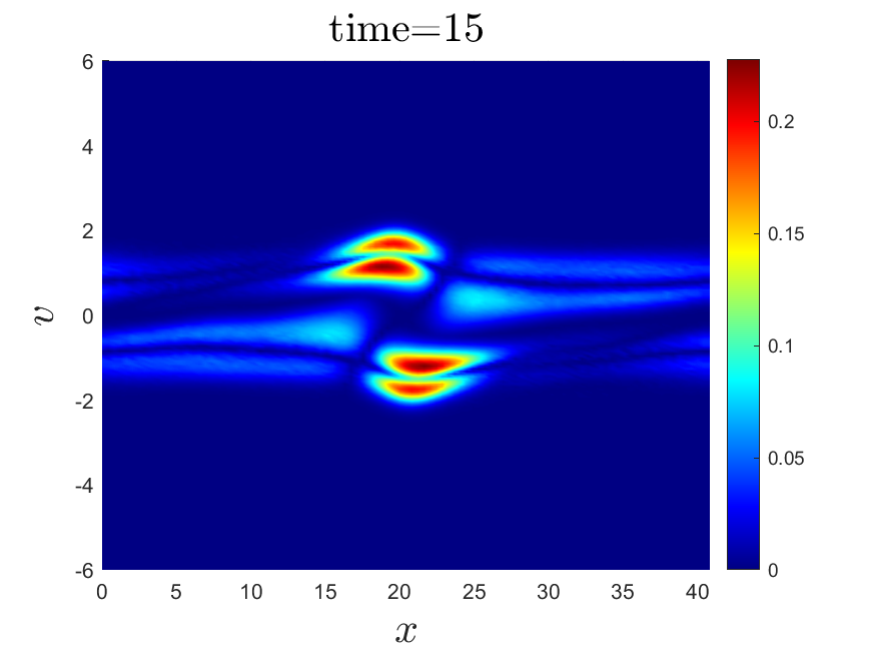}
	\includegraphics[width = 0.325\linewidth]{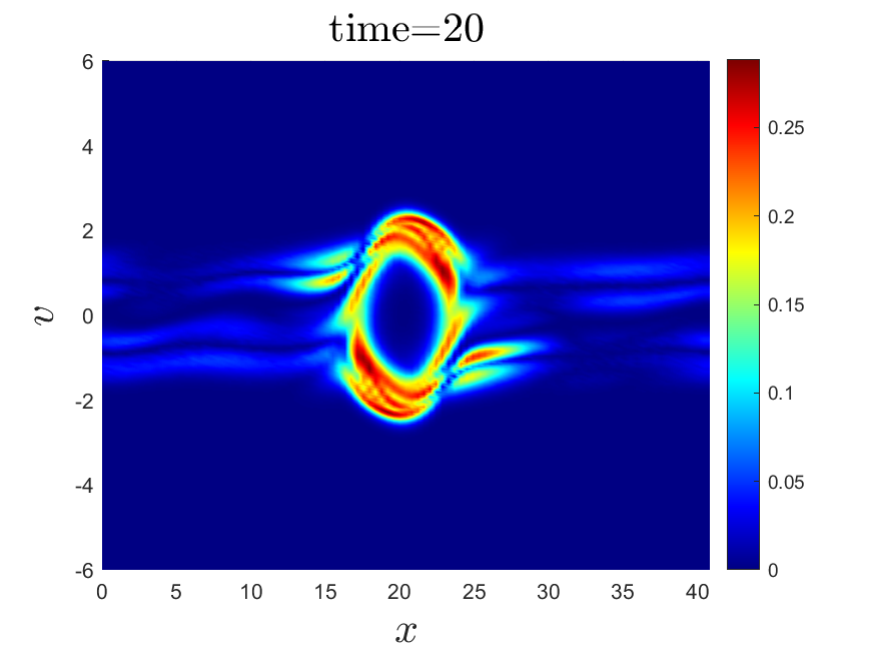}
	\caption{\small{\textbf{Test 3: nonlinear two stream instability}. We report $\mathbb{E}_{\z}[f(x,v,t,\z)]$  (top row) and $\textrm{Var}_{\z}[f(x,v,t,\z)]$ (bottom row) at fixed times $t=0,15,20$, in collisionless nonlinear two stream instability. We choose $k=2/13$, $N=5 \cdot 10^7$ particles, $M=5$ and $\Delta t=0.1$. We considered the initial condition \eqref{eq:init_nonlinear_twostream} with $\alpha(\z) = 4\cdot 10^{-2} + 2\cdot 10^{-2}\z$, $\z\sim\mathcal{U}([0,1])$.}}
	\label{fig:test_3_non_linear_two_streams_epsi_inf}
\end{figure}

In the collisionless scenario, the nonlinear phenomenon of the electron trapping causes stronger distortions of the distribution $f(x,v,t,\z)$, which can not be treated analytically, as we can observe from $\mathbb{E}_{\z}[f(x,v,t,\z)]$ in Figure \ref{fig:test_3_non_linear_two_streams_epsi_inf} at times $t = 0,15,20$. In the same figure we report the values of $\textrm{Var}_{\z}[f(x,v,t,\z)]$ computed at the same times from which we get information on the deviation from expected trends. 

\begin{figure}
	\centering
	\includegraphics[width = 0.325\linewidth]{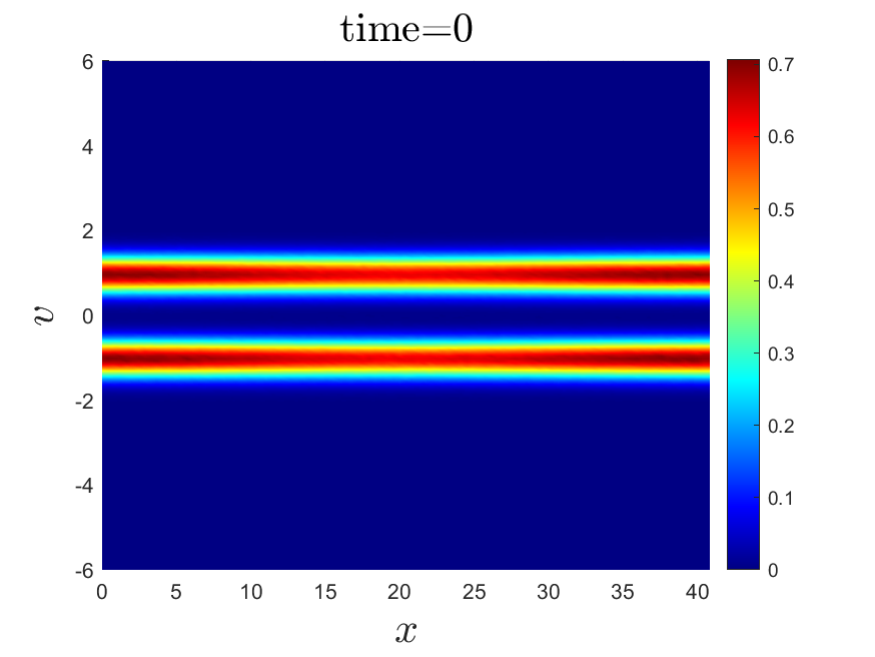}
	\includegraphics[width = 0.325\linewidth]{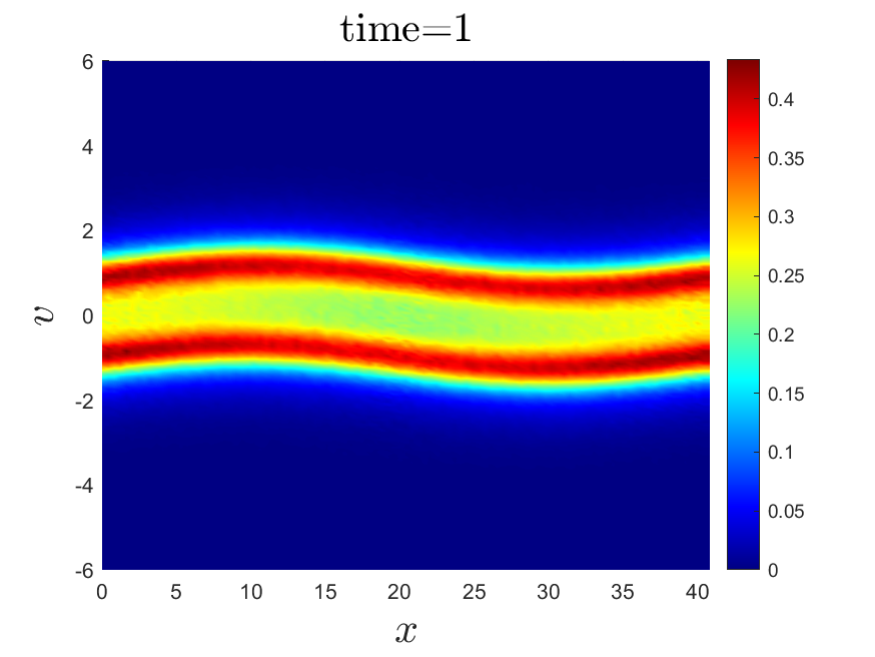}
	\includegraphics[width = 0.325\linewidth]{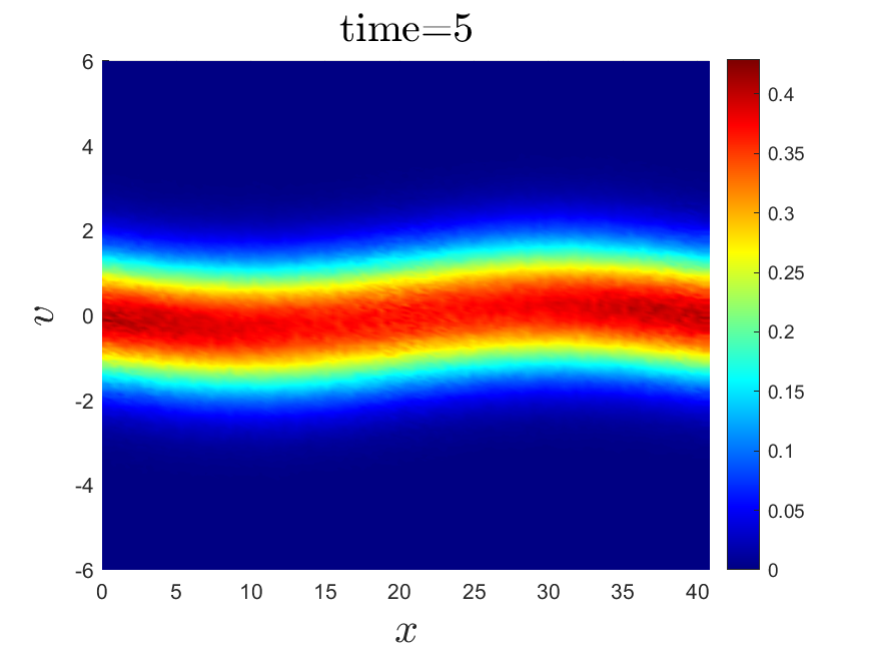}
	\includegraphics[width = 0.325\linewidth]{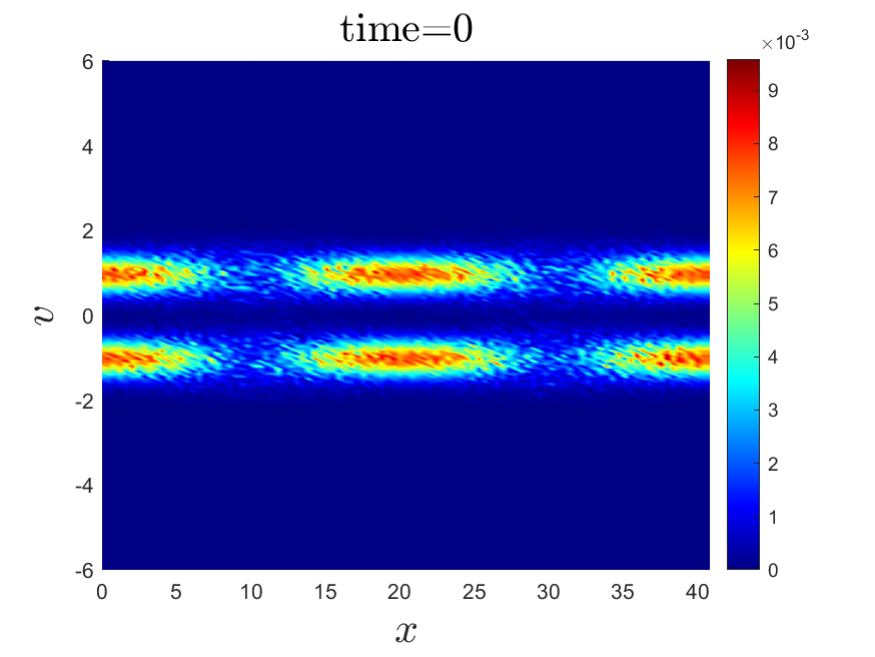}
	\includegraphics[width = 0.325\linewidth]{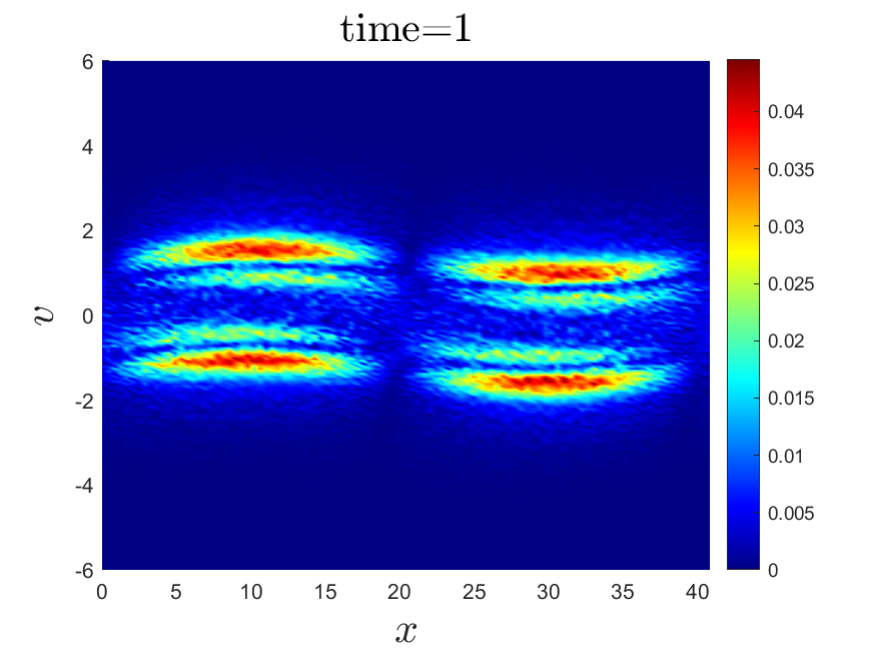}
	\includegraphics[width = 0.325\linewidth]{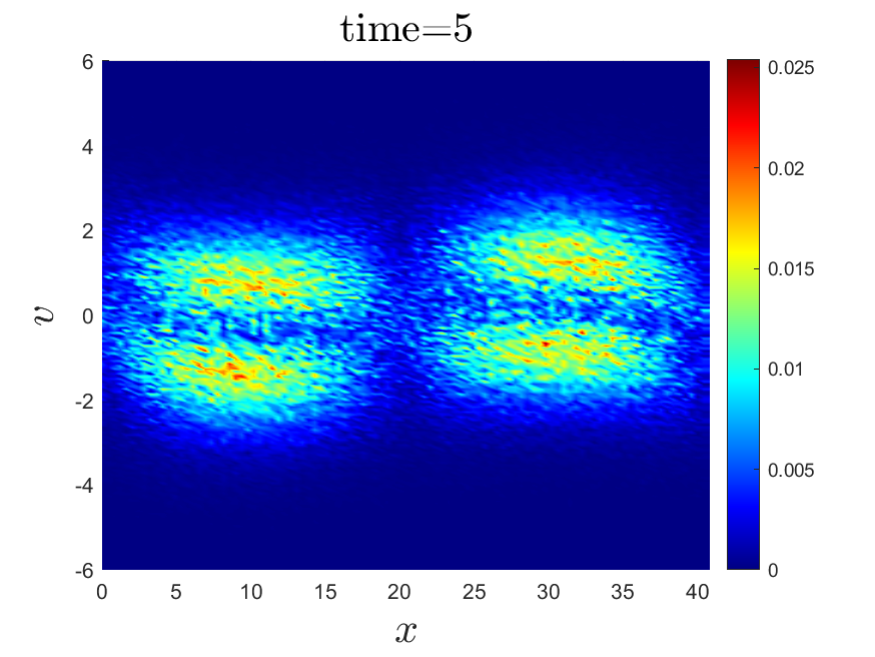}
	\caption{\small{\textbf{Test 3: nonlinear two stream instability}. We report $\mathbb{E}_{\z}[f(x,v,t,\z)]$  (top row) and $\textrm{Var}_{\z}[f(x,v,t,\z)]$ (bottom row) at fixed times $t=0,1,5$, in nonlinear two stream instability with $\nu=1$. We choose $k=2/13$, $N=5 \cdot 10^7$ particles, $M=5$ and $\Delta t=0.1$. We considered the initial condition \eqref{eq:init_nonlinear_twostream} with $\alpha(\z) = 4\cdot 10^{-2} + 2\cdot 10^{-2}\z$, $\z\sim\mathcal{U}([0,1])$. }}
	\label{fig:test_3_non_linear_two_streams_epsi_1}
\end{figure}

\begin{figure}
	\centering
	\includegraphics[width = 0.325\linewidth]{Immagini/non_lin_two_stream_epsi_1_ex_F0}
	\includegraphics[width = 0.325\linewidth]{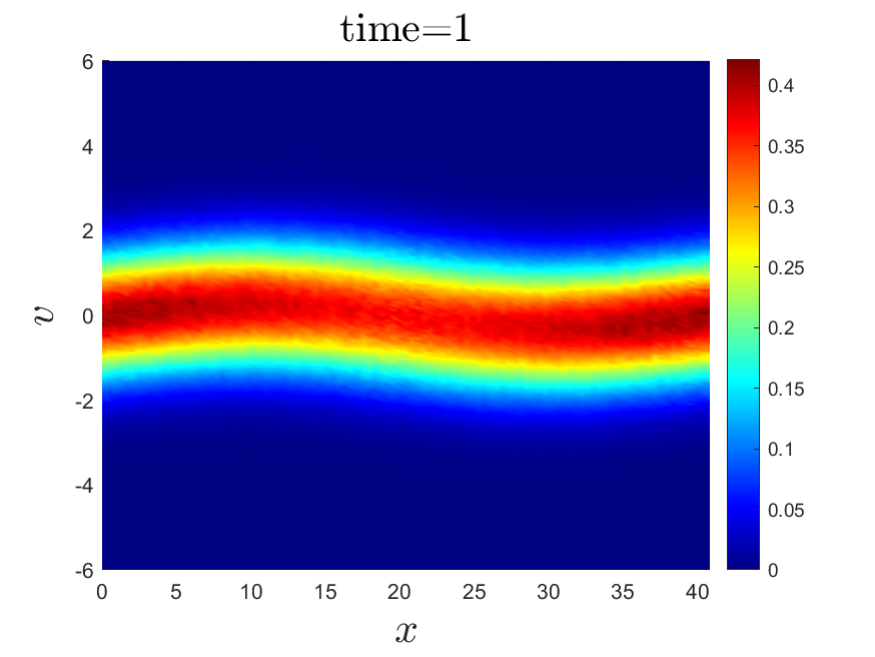}
	\includegraphics[width = 0.325\linewidth]{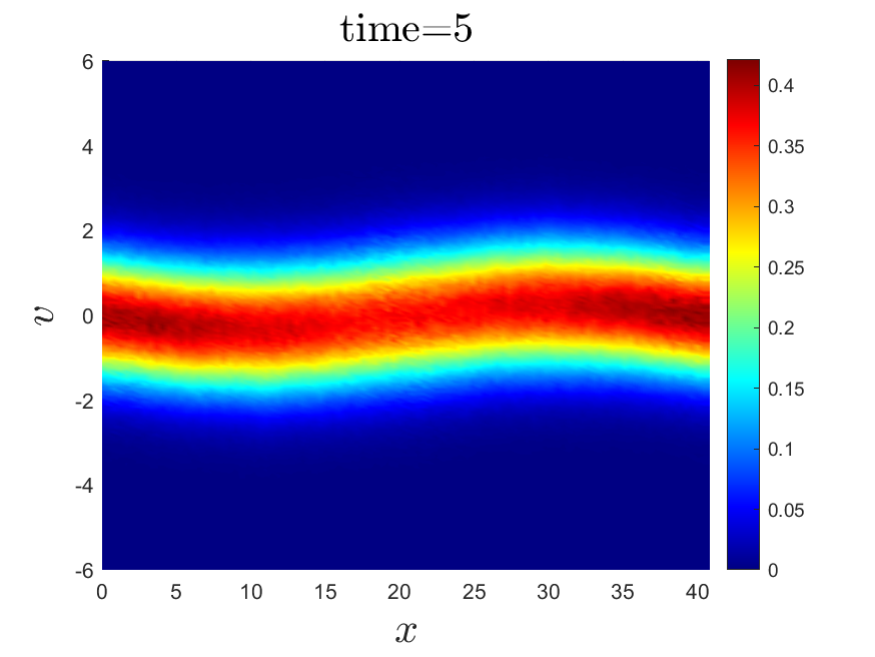}
	\includegraphics[width = 0.325\linewidth]{Immagini/non_lin_two_stream_epsi_1_var_F0}
	\includegraphics[width = 0.325\linewidth]{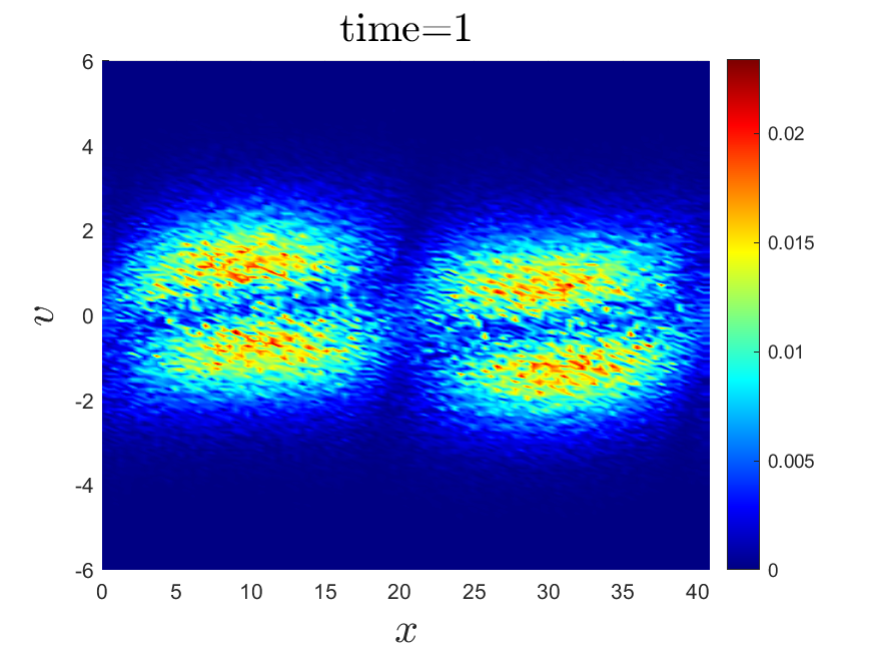}
	\includegraphics[width = 0.325\linewidth]{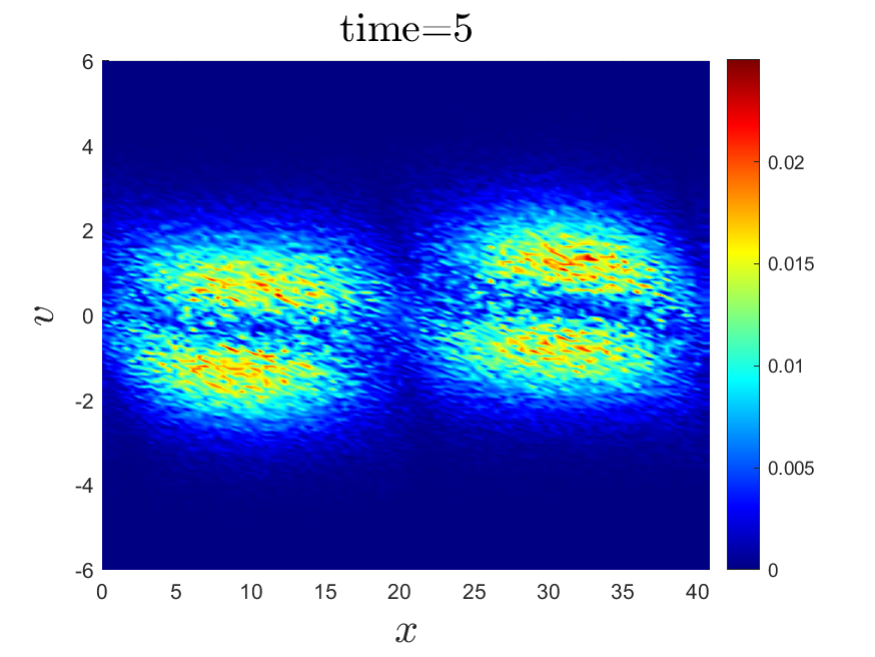}
	\caption{\small{\textbf{Test 3: nonlinear two stream instability}. We report $\mathbb{E}_{\z}[f(x,v,t,\z)]$  (top row) and $\textrm{Var}_{\z}[f(x,v,t,\z)]$ (bottom row) at fixed times $t=0,1,5$, in nonlinear two stream instability with $\nu=10^3$. We choose $k=2/13$, $N=5 \cdot 10^7$ particles, $M=5$ and $\Delta t=0.1$. We considered the initial condition \eqref{eq:init_nonlinear_twostream} with $\alpha(\z) = 4\cdot 10^{-2} + 2\cdot 10^{-2}\z$, $\z\sim\mathcal{U}([0,1])$.}}
	\label{fig:test_3_non_linear_two_streams_epsi_0_001}
\end{figure}

If we consider a non-zero collision frequency $\nu$, the interactions among the particles balance the electron trapping phenomenon letting the system reach the Maxwellian equilibrium quickly, as we can see from Figure \ref{fig:test_3_non_linear_two_streams_epsi_1} and \ref{fig:test_3_non_linear_two_streams_epsi_0_001}.

\subsection{Test 4: Sod shock tube}
To evaluate the capability of the sG particle scheme to capture the limit Euler-Poisson system \eqref{eq:euler_poisson}, we perform the Sod shock tube test^^>\cite{sod1978} for different collisional regimes, {comparing our results with a collocation method applied to a Lax-Fridrichs and a WENO scheme \cite{shu2009}}. The test is a typical Riemann problem, namely it consists in a piecewise constant initial value problem with a discontinuity in the space domain. In particular, we consider two different settings: the first with an uncertain initial temperature, the second with an uncertain position of the interface.

\paragraph{Uncertain initial temperature.}
We consider the initial distribution
\[
f_0(x,v,\z) = \dfrac{\rho_0(x)}{\sqrt{2\pi T_0(x,\z)}} e^{-\dfrac{v^2}{2T_0(x,\z)}},
\]
with $x\in[0,1]$. The mass $\rho_0(x)$ and the uncertain temperature $T_0(x,\z)$ are initialized as 
\be\label{eq:unc_temp}
\begin{split}
	\rho_0(x) = 1, \qquad T_0(x,\z)=1+\alpha(\z) \qquad &\textrm{if}\quad 0<x<0.5  \\
	\rho_0(x) = 0.125, \qquad T_0(x,\z)=0.8+\alpha(\z) \qquad &\textrm{if}\quad 0.5<x<1
\end{split}
\ee
with $\alpha(\z) = 0.25 \z$, $\z\sim\mathcal{U}([0,1])$. We choose $N=10^7$, $M=5$, $\Delta t=0.01$ and we consider different collisional regimes corresponding to the choices $\nu=1, 10^3$. We implement Dirichelet boundary conditions for the Poisson equation and reflecting boundary conditions for the particles according to \eqref{eq:rbc_x} and \eqref{eq:rbc_v}. 

\begin{figure}
	\centering
	\includegraphics[width = 0.4\linewidth]{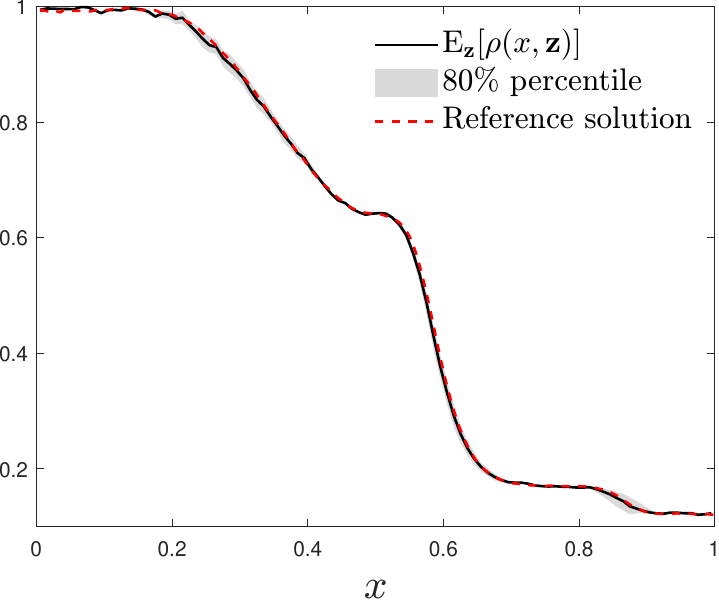}
	\includegraphics[width = 0.4\linewidth]{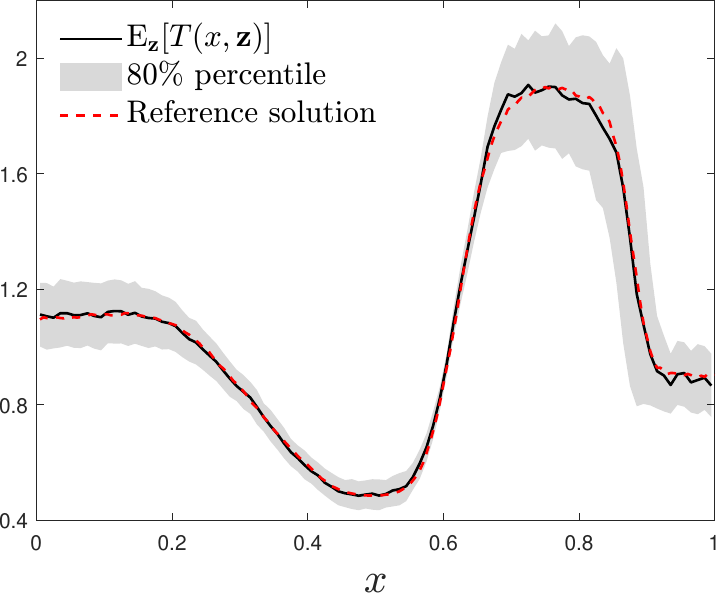}
	\includegraphics[width = 0.4\linewidth]{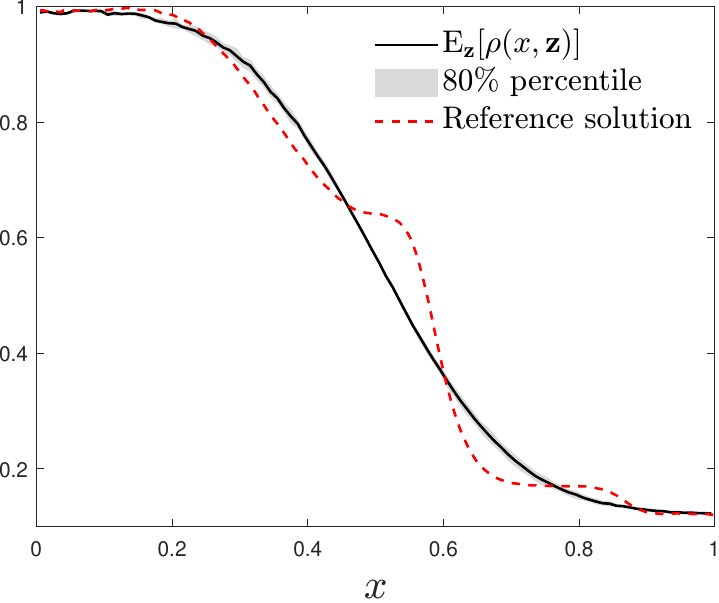}
	\includegraphics[width = 0.4\linewidth]{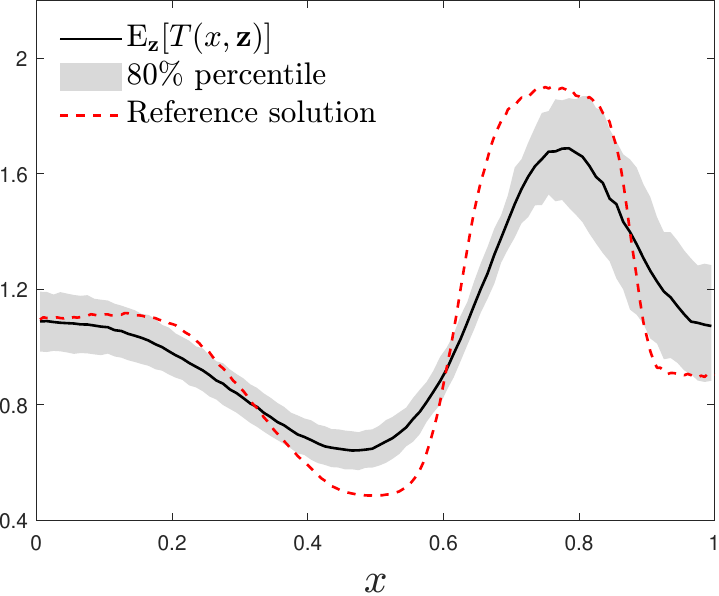}
	\caption{\small{\textbf{Test 4: Sod shock tube}. We report $\mathbb{E}_{\z}[\rho(x,t,\z)]$ (left) and $\mathbb{E}_{\z}[T(x,t,\z)]$ (right) for $\nu=10^3$ (top) and $\nu=1$ (bottom), at fixed times $t=0.15$, for the Sod shock test with uncertain initial temperature. We choose $N=10^7$, $M=5$, $\Delta t=0.01$. We considered the initial condition \eqref{eq:unc_temp} with $\alpha(\z) = 0.25 \z$, $\z\sim\mathcal{U}([0,1])$. Reference solution computed with $N=10^8$ particles.}}
	\label{fig:test_4_hydro}
\end{figure}

We compute a reference solution with $N=10^8$ particles, $N_\ell=200$ cells, $\nu=10^4$ and the other parameters as before, to compare the results of the two collisional regimes. In particular, we compare the expectations of the mass $\rho(x,\z)$ and temperature $T(x,\z)$ at the final time $t=0.15$, with respect to the reference solution. In Figure \ref{fig:test_4_hydro} we observe that for $\nu=10^3$ (top row) the mean values together with the confidence intervals fit the reference solution, while for $\nu=1$ (bottom row) the frequency of collisions is obviously not enough to justify the fluid limit.

\begin{figure}
	\centering
	\includegraphics[width = 0.4\linewidth]{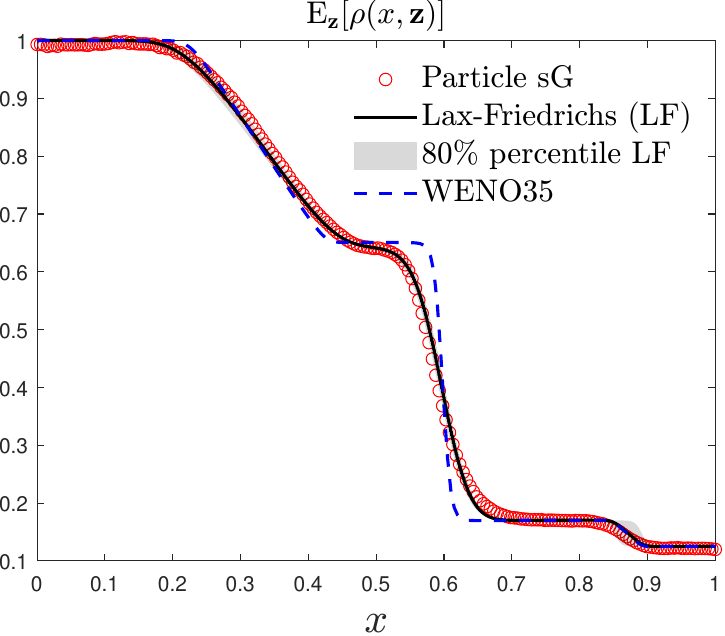}
	\includegraphics[width = 0.4\linewidth]{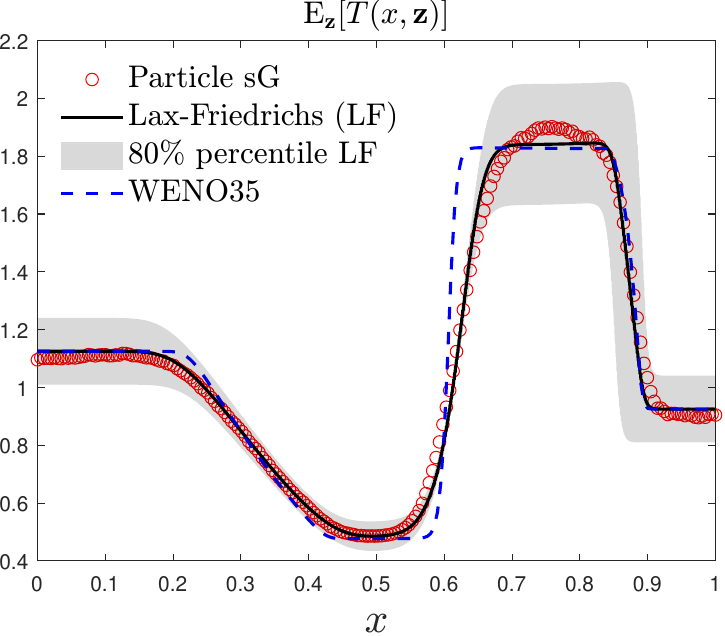}
	\caption{\small{\textbf{Test 4: Sod shock tube}. Comparison of the particle sG reference solution with Lax-Friedrichs and WENO schemes. We considered the initial condition \eqref{eq:unc_temp} with $\alpha(\z) = 0.25 \z$, $\z\sim\mathcal{U}([0,1])$. The particle solution is computed with $N=10^8$, $M=5$ and $\Delta t=0.01$. Lax-Friedrichs is solved with $1500$ cells and CFL number $0.1$, WENO with $200$ cells and CFL number $0.5$. In both cases we use a stochastic collocation approach with $11$ nodes.}}
	\label{fig:test_4_hydro_comp}
\end{figure}

{
We solve the same uncertain Riemann problem with a first order in time and space Lax-Friedrichs (LF) method, and a third order in time, fifth order in space WENO scheme. For the uncertain parameter, we adopted a stochastic collocation method with $11$ collocation nodes, LF is solved with $1500$ cells and a CFL number equal to $0.1$, WENO35 with $200$ cells and CFL number $0.5$. The results are summarized in Figure \ref{fig:test_4_hydro_comp}. We observe that the reference particle sG solution (red circles) is in good accordance with the Lax-Friedrichs results (solid black line), since they both are first order.
}

\paragraph{Uncertain interface position.}
We consider now the distribution at the initial time 
\[
f_0(x,v,\z) = \dfrac{\rho_0(x,\z)}{\sqrt{2\pi T_0(x,\z)}} e^{-\dfrac{v^2}{2T_0(x,\z)}},
\]
with $x\in[0,1]$. In this scenario, we design an uncertain initial interface, i.e.
\be\label{eq:unc_int}
\begin{split}
	\rho_0(x,\z) = 1, \qquad T_0(x,\z)=1 \qquad &\textrm{if}\quad 0<x<0.5 +\alpha(\z)  \\
	\rho_0(x,\z) = 0.125, \qquad T_0(x,\z)=0.8 \qquad &\textrm{if}\quad 0.5+\alpha(\z)<x<1
\end{split}
\ee
with $\alpha(\z) = -0.05+0.1\z$, $\z\sim\mathcal{U}([0,1])$. Again, we choose $N=10^7$, $M=5$, $\Delta t=0.01$ and we consider different collisional regimes corresponding to the choices $\nu=1, 10^3$. The boundary conditions for the particles and the Poisson equation are the same as the previous test.

\begin{figure}
	\centering
	\includegraphics[width = 0.4\linewidth]{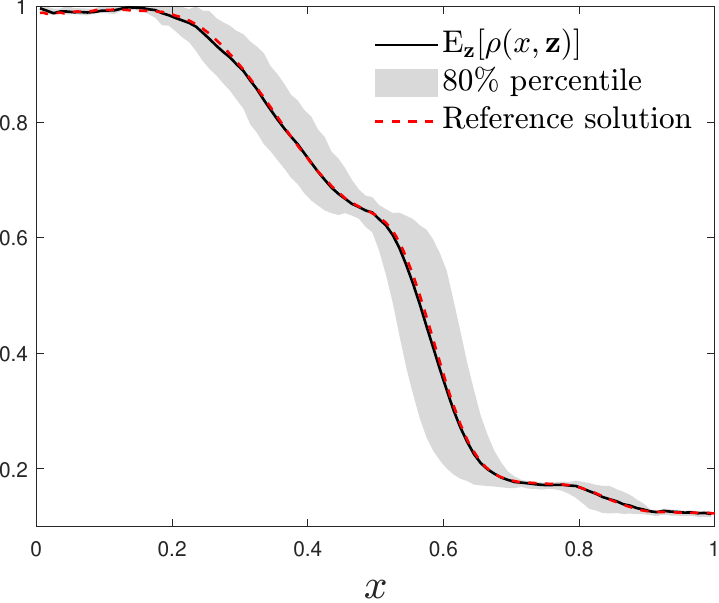}
	\includegraphics[width = 0.4\linewidth]{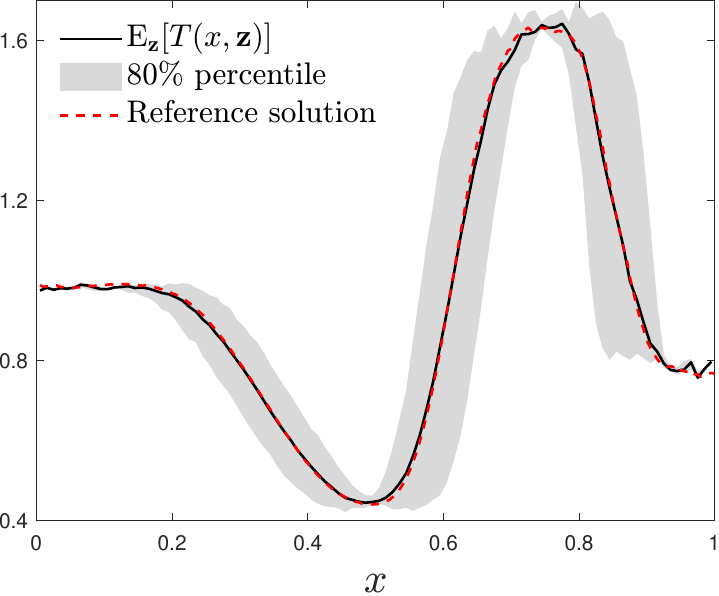}
	\includegraphics[width = 0.4\linewidth]{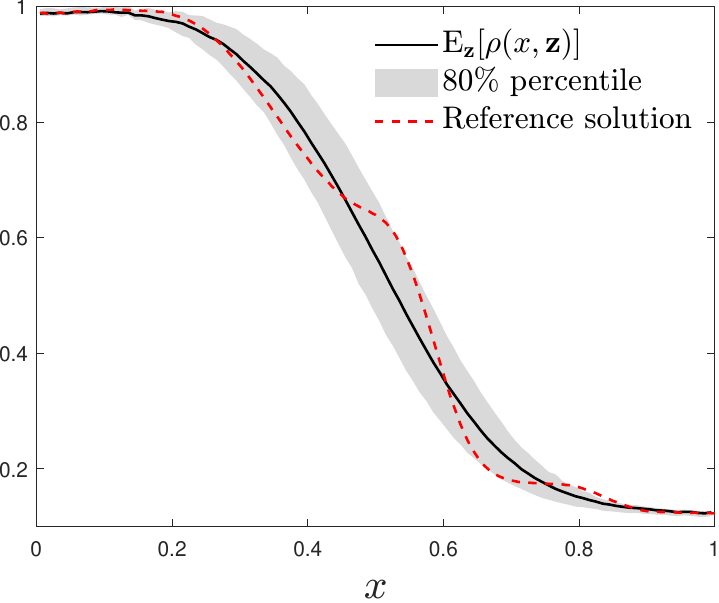}
	\includegraphics[width = 0.4\linewidth]{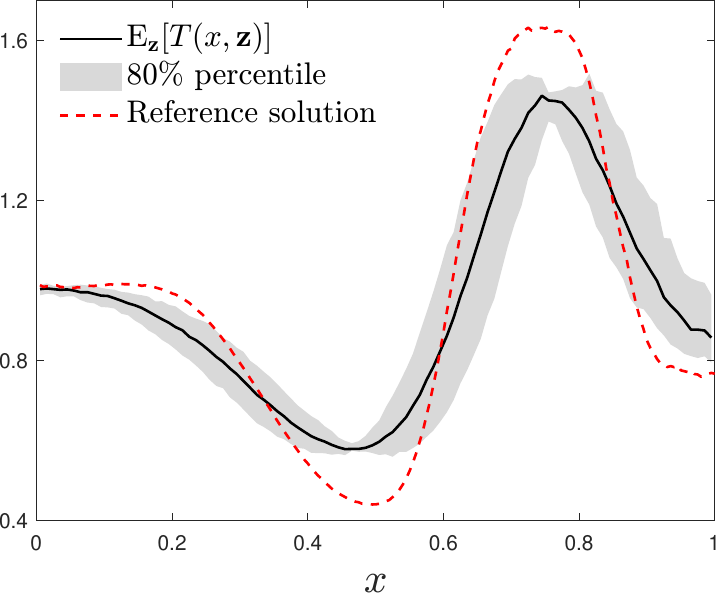}
	\caption{\small{\textbf{Test 4: Sod shock tube}. We report $\mathbb{E}_{\z}[\rho(x,t,\z)]$ (left) and $\mathbb{E}_{\z}[T(x,t,\z)]$ (right) for $\nu=10^3$ (top) and $\nu=1$ (bottom), at fixed times $t=0.15$, for the Sod shock test with uncertain interface position. We choose $N=10^7$, $M=5$, $\Delta t=0.01$. We considered the initial condition \eqref{eq:unc_int} with $\alpha(\z) = -0.05 + 0.1\z$, $\z\sim\mathcal{U}([0,1])$. Reference solution computed with $N=10^8$ particles.}}
	\label{fig:test_4_hydro2}
\end{figure}

The reference solution is computed with $N=10^8$ particles, $N_\ell=200$ cells, $\nu=10^4$ and the other parameters as before. In Figure \ref{fig:test_4_hydro2} we show the expectations and the confidence intervals of the mass and temperature at the final time $t=0.15$. We observe that in the collisional regime, assumed here to be characterized by the frequency $\nu=10^3$ (top row), the results are in agreement with the reference solution. For $\nu=1$ (bottom row) clearly the collision frequency is far from the fluid limit and the profiles computed from the kinetic density differ from the ones of reference solution. 

\begin{figure}
	\centering
	\includegraphics[width = 0.4\linewidth]{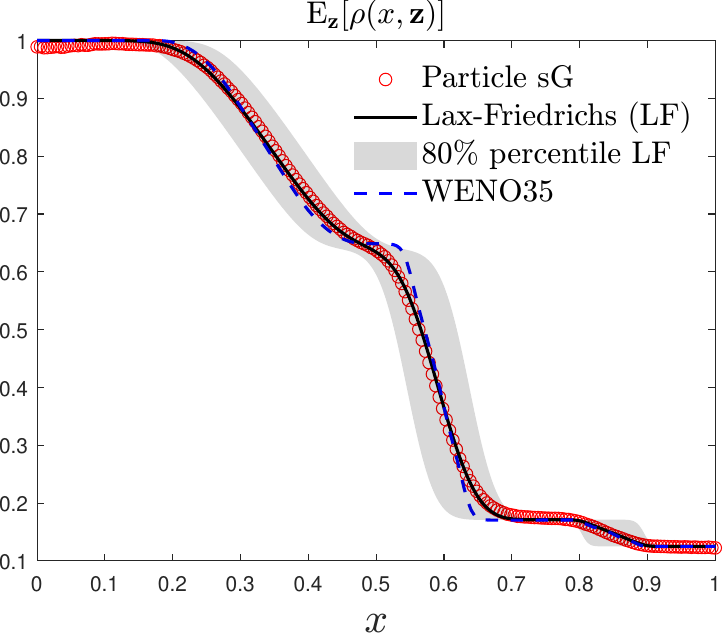}
	\includegraphics[width = 0.4\linewidth]{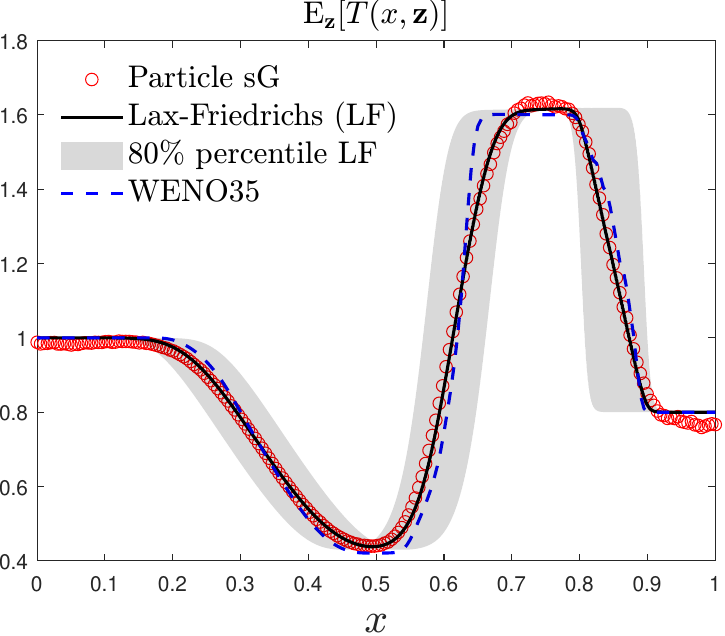}
	\caption{\small{\textbf{Test 4: Sod shock tube}. Comparison of the particle sG reference solution with Lax-Friedrichs and WENO schemes. We considered the initial condition \eqref{eq:unc_int} with $\alpha(\z) = -0.05 + 0.1\z$, $\z\sim\mathcal{U}([0,1])$. The particle solution is computed with $N=10^8$, $M=5$ and $\Delta t=0.01$. Lax-Friedrichs is solved with $1500$ cells and CFL number $0.1$, WENO with $200$ cells and CFL number $0.5$. In both cases we use a stochastic collocation approach with $11$ nodes.}}
	\label{fig:test_4_hydro_comp2}
\end{figure}

{
As in the previous case, we compare the particle sG solver with a first order in time and space Lax-Friedrichs (LF) method, and a third order in time, fifth order in space WENO scheme. The computational setting is the same as before: we adopted a stochastic collocation method with $11$ collocation nodes, LF is solved with $1500$ cells and a CFL number equal to $0.1$, WENO35 with $200$ cells and CFL number $0.5$. In Figure \ref{fig:test_4_hydro_comp2} we may observe that the reference particle sG solution (red circles) is in good accordance with the Lax-Friedrichs results (solid black line).
}


\section{Conclusions}
We presented a novel numerical method based on a sG projection of the corresponding particle dynamics for uncertainty quantification in collisional plasma physics. From a mathematical viewpoint, the system under investigation is described by a Vlasov-Poisson kinetic model with BGK-type interactions. The numerical strategy, unlike sG approximations based on a direct discretization of the system of PDEs, preserves the main physical properties of the system, including the positivity of the solution. In addition, thanks to a splitting approach and an exact representation of the collisional process, the method satisfies the asymptotic-preserving property. As a consequence, in the fluid-dynamic limit, it originates a sG particle method for the corresponding Euler-Poisson system. A collection of numerical examples, ranging from linear and nonlinear Landau damping to Sod shock tube, highlighted the spectral accuracy for smooth solutions in the random space and the validity of the present approach for problems with different space and time scales in presence of uncertainty. Future research directions will be devoted to the introduction of a Landau type collision term for which the application of the sG particle method is not straightforward.

\section*{Acknowledgments}
This work has been written within the
activities of GNCS and GNFM groups of INdAM (National Institute of
High Mathematics). L.P. acknowledges the partial support of MIUR-PRIN Project 2017, No. 2017KKJP4X “Innovative numerical methods for evolutionary partial differential equations and applications”. M.Z. acknowledges partial support of MUR-PRIN2020 Project, No. 2020JLWP23 “Integrated mathematical approaches to socio-epidemiological dynamics”. The authors acknowledge the support of the Banff International Research Station (BIRS) for the Focused Research Group [22frg198] “Novel perspectives in kinetic equations for emerging phenomena”, July 17-24, 2022, where part of this work was done.

\appendix

\section{{Initial sampling}} \label{appendix}
{
Given the initial condition $f_0(x,v,\z)$, we consider the empirical distribution
\[
f^0_N(x,v,\z)=\frac{m(\z)}{N} \sum_{i=1}^{N}  S( x - x^0_i(\z)  ) \otimes S( v - v^0_i(\z)  ),
\]
in a way that formally $f^0_N(x,v,\z)\to f_0(x,v,\z)$ as $N\to\infty$. However, as already pointed out in \cite{Pareschi2020}, a direct application of standard sampling techniques may be nontrivial since we have to perform a gPC projection to obtain the coefficients as presented in Section \ref{subsec:sG}. }

\textcolor{blue}{
\begin{figure}[htb]
	\centering
	\begin{minipage}{.9\linewidth}
		\begin{algorithm}[H]  
			\footnotesize
			\caption{\small{sG projection of the initial positions} }
			\begin{itemize}
				\item Generate a set of Gaussian nodes $\{z_k\}$, for $k=1,\dots,K$, from the distribution of the parameters $p(\z)$;
				\item discretize the spatial domain $[L_{-},L_{+}]$ in $N_\ell$ equally spaced cells in a way that $x_0=L_{-},\dots,x_{\ell},\dots,x_{N_\ell}=L_{+}$;
				\item for $k=1$ to $K$:
				\begin{itemize}
					\item compute the mass fraction of every cell $\rho^0_{\ell, k} = \rho^0_\ell(z_k)/m(z_k)$, with $m(z_k)=\int \rho^0_\ell(z_k) dx$;
					\item set $N_{\ell, k}=N\rho^0_{\ell, k}$, number of particles to sample in cell;
					\item generate $N_{\ell, k}$ uniformly distributed positions $\{x^{0,M}_i(z_k)\}$ in the interval $[x_{\ell-1}, x_{\ell}]$, for every $\ell=1,\dots,N_\ell$;						
				\end{itemize}
				end for;
				\item for $h=0,\dots,M$:
				\begin{itemize}
					\item compute the projections $\hat{x}^{0}_{i,h}$ according to \eqref{eq:proj};
				\end{itemize}
				end for.
			\end{itemize}
		\end{algorithm}
	\end{minipage}
\end{figure}}
\textcolor{blue}{
\begin{figure}[htb]
	\centering
	\begin{minipage}{.9\linewidth}
		\begin{algorithm}[H]  
			\footnotesize
			\caption{\small{sG projection of the initial velocities} }
			\begin{itemize}
				\item Generate a set of Gaussian nodes $\{z_k\}$, for $k=1,\dots,K$, from the distribution of the parameters $p(\z)$;
				\item for $i=1$ to $N$:
				\begin{itemize}
					\item generate a reference sample of $\tilde  v^{0}_i$ using suitable sampling method from the expected marginal $\mathbb E_{\z}[f_v]$ of $f_0(x,v,\z)$ in $x$;
					\item compute $v^{0,M}_i(z_k) = \sqrt{\frac{T(z_k)}{\tilde T}}\tilde  v^{0}_i$ for each $k=1,\dots,K$, being $\tilde  T$ the temperature of $\mathbb E_{\z}[f_v]$					
				\end{itemize}
				end for;
				\item for $h=0,\dots,M$:
				\begin{itemize}
					\item compute the projections $\hat{v}^{0}_{i,h}$ according to \eqref{eq:proj};
				\end{itemize}
				end for.
			\end{itemize}
		\end{algorithm}
	\end{minipage}
\end{figure}}

{
Indeed, we first need to construct a random sample $x^{0,M}_i(\z)$, $v^{0,M}_i(\z)$, for $i=1,\dots,N$ such that
\[
f^0_N(x,v,\z) \approx \frac{m(\z)}{N} \sum_{i=1}^{N}  S( x - x^{0,M}_i(\z)  ) \otimes S( v - v^{0,M}_i(\z)  ),
\]
where $x^{0,M}_i(\z)$ and $v^{0,M}_i(\z)$ are given by \eqref{eq:gPCE}. Hence, by direct integration we can compute the coefficients of the polynomial chaos expansion according to \eqref{eq:proj}.}

{
Since the integrals in the random space are approximated through a Gaussian quadrature rule, we want to generate a sample of positions and velocities in every Gauss collocation node $\{z_k\}_{k=1}^{K}$, in a way that
\[
\begin{split}
	&\hat{x}^0_{i,h} = \int_{\Omega} x^0_i(\z) \Psi_h(\z) p(\z) d\z \approx \sum_{k=1}^{M} x^{0,M}_i(z_k) w_k \Psi_h(z_k) , \\
	&\hat{v}^0_{i,h} = \int_{\Omega} v^0_i(\z) \Psi_h(\z) p(\z) d\z \approx \sum_{k=1}^{M} v^{0,M}_i(z_k) w_k \Psi_h(z_k).
\end{split}
\]
As a consequence, given $f_0(x,v,\z)$ and $p(\z)$, distribution of the random parameters, we desire a set of samples
\[
\begin{split}
& X^{1} = \{x^{0,M}_1(z_1), \dots, x^{0,M}_N(z_1)\},\; \dots,\; X^{K} = \{x^{0,M}_1(z_K), \dots, x^{0,M}_N(z_K)\} \\
& V^{1} = \{v^{0,M}_1(z_1), \dots, v^{0,M}_N(z_1)\},\; \dots,\; V^{K} = \{v^{0,M}_1(z_K), \dots, v^{0,M}_N(z_K)\},
\end{split}
\]
such that $\{x^{0,M}_i(z_k),\,v^{0,M}_i(z_k)\}$, for $k=1,\dots,K$, represent the different realizations of the same $i$-th particle $\{x^{0,M}_i(\cdot),\,v^{0,M}_i(\cdot)\}$. We observe also that we can sample independently positions and velocities since the initial distribution functions we considered are tensorised. }

{
At the numerical level, one possible strategy in one-dimension is to sample $X^{k}$ and $V^{k}$ independently with respect to the nodes $z_k$, with any sampling method, and then to sort positions and velocities for any $k=1,\dots,K$. In higher dimensions, the ordering technique is not possible and we need to solve an assignment problem in order to link particles. In this work, we focused on the 1D case, and we adopted the following techniques.}
{
To observe the physical phenomena presented in Section \ref{sec:test} we need to sample the density $\rho_0(\z)$ carefully by reducing statistical fluctuations. To this aim, we adopted a stratified sampling method in each Gaussian node by sorting the vectors $X^k$ for every $k$ at the end.}

{
Note that, if we consider continuous probability distributions with uncertainties in the variance, e.g. a Maxwellian with uncertain temperature, we can overcome the sorting method thanks to the scaling property of the second order moment (see \cite{Pareschi2020}).}

\bibliographystyle{abbrv}
\bibliography{UQ_Plasmi.bib}

\begin{thebibliography}{10}

\bibitem{Andries2002}
P.~Andries, K.~Aoki, and B.~Perthame.
\newblock A consistent {BGK}-type model for gas mixtures.
\newblock {\em J. Stat. Phys.}, 106:993--1018, 2002.

\bibitem{BGK}
P.~L. Bhatnagar, E.~P. Gross, and M.~Krook.
\newblock A {M}odel for {C}ollision {P}rocesses in {G}ases. {I}. {S}mall
  {A}mplitude {P}rocesses in {C}harged and {N}eutral {O}ne-{C}omponent
  {S}ystems.
\newblock {\em Phys. Rev.}, 94(3):511, 1954.

\bibitem{Dimarco2008}
R.~Caflisch, C.~Wang, G.~Dimarco, B.~Cohen, and A.~Dimits.
\newblock A hybrid method for accelerated simulation of {C}oulomb collisions in
  a plasma.
\newblock {\em Multiscale Model. Simul.}, 7(2):865--887, 2008.

\bibitem{Campospinto2014}
M.~Campos-Pinto, E.~Sonnendr{\"u}cker, A.~Friedman, D.~P. Grote, and S.~M.
  Lund.
\newblock Noiseless {V}lasov--{P}oisson simulations with linearly transformed
  particles.
\newblock {\em J. Comput. Phys.}, 275:236--256, 2014.

\bibitem{Jingwei2020}
J.~A. Carrillo, J.~Hu, L.~Wang, and J.~Wu.
\newblock A particle method for the homogeneous {L}andau equation.
\newblock {\em J. Comput. Phys. X}, 7:100066, 24, 2020.

\bibitem{Carrillo2022}
J.~A. Carrillo, S.~Jin, and Y.~Tang.
\newblock Random batch particle methods for the homogeneous {L}andau equation.
\newblock {\em Commun. Comput. Phys.}, 31(4):997--1019, 2022.

\bibitem{CPZ19}
J.~A. Carrillo, L.~Pareschi, and M.~Zanella.
\newblock Particle based g{PC} methods for mean-field models of swarming with
  uncertainty.
\newblock {\em Commun. Comput. Phys.}, 25(2):508--531, 2019.

\bibitem{CZ19}
J.~A. Carrillo and M.~Zanella.
\newblock Monte {C}arlo g{PC} methods for diffusive kinetic flocking models
  with uncertainties.
\newblock {\em Vietnam J. Math.}, 47(4):931--954, 2019.

\bibitem{chacon2016}
E.~Chacon-Golcher, S.~A. Hirstoaga, and M.~Lutz.
\newblock Optimization of {P}article-{I}n-{C}ell simulations for
  {V}lasov--{P}oisson system with strong magnetic field.
\newblock {\em ESAIM: Proc.}, 53:177--190, 2016.

\bibitem{Chen1974}
F.~F. Chen.
\newblock {\em Introduction to {P}lasma {P}hysics and {c}ontrolled {F}usion}.
\newblock Plenum {P}ress, New {Y}ork and {L}ondon, 2nd edition, 1974.

\bibitem{Chung2020}
S.~W. Chung, S.~D. Bond, E.~C. Cyr, and J.~B. Freund.
\newblock Regular sensitivity computation avoiding chaotic effects in
  particle-in-cell plasma methods.
\newblock {\em J. Comput. Phys.}, 400:108969, 32, 2020.

\bibitem{crestetto2012}
A.~Crestetto, N.~Crouseilles, and M.~Lemou.
\newblock Kinetic/fluid micro-macro numerical schemes for
  {V}lasov-{P}oisson-{BGK} equation using particles.
\newblock {\em Kinet. Relat. Models}, 5(4):787--816, 2012.

\bibitem{Crouseilles2004}
N.~Crouseilles and F.~Filbet.
\newblock Numerical approximation of collisional plasmas by high order methods.
\newblock {\em J. Comput. Phys.}, 201(2):546--572, 2004.

\bibitem{Dai2022}
D.~Dai, Y.~Epshteyn, and A.~Narayan.
\newblock Hyperbolicity-preserving and well-balanced stochastic {G}alerkin
  method for two-dimensional shallow water equations.
\newblock {\em J. Comput. Phys.}, 452:110901, 2022.

\bibitem{despres2013}
B.~Despr{\'e}s, G.~Po{\"e}tte, and D.~Lucor.
\newblock Robust uncertainty propagation in systems of conservation laws with
  the entropy closure method.
\newblock In {\em Uncertainty {Q}uantification in {C}omputational {F}luid
  {D}ynamics}, pages 105--149. Springer, 2013.

\bibitem{Dimarco2010}
G.~Dimarco, R.~Caflisch, and L.~Pareschi.
\newblock Direct simulation {M}onte {C}arlo schemes for {C}oulomb interactions
  in plasmas.
\newblock {\em Commun. Appl. Ind. Math.}, 1(1):72--91, 2010.

\bibitem{Dimarco2015}
G.~Dimarco, Q.~Li, L.~Pareschi, and B.~Yan.
\newblock Numerical methods for plasma physics in collisional regimes.
\newblock {\em J. Plasma Phys.}, 81(1):1--31, 2015.

\bibitem{Dimarco2014}
G.~Dimarco, L.~Mieussens, and V.~Rispoli.
\newblock An asymptotic preserving automatic domain decomposition method for
  the {V}lasov-{P}oisson-{BGK} system with applications to plasmas.
\newblock {\em J. Comput. Phys.}, 274:122--139, 2014.

\bibitem{PD2014}
G.~Dimarco and L.~Pareschi.
\newblock Numerical methods for kinetic equations.
\newblock {\em Acta Numer.}, 23:369--520, 2014.

\bibitem{DP2020}
G.~Dimarco and L.~Pareschi.
\newblock Multiscale variance reduction methods based on multiple control
  variates for kinetic equations with uncertainties.
\newblock {\em Multiscale Model. Simul.}, 18(1):351--382, 2020.

\bibitem{Jin2019}
Z.~Ding and S.~Jin.
\newblock Random regularity of a nonlinear {L}andau damping solution for the
  {V}lasov-{P}oisson equations with random inputs.
\newblock {\em Int. J. Uncertain. Quantif.}, 9(2):123--142, 2019.

\bibitem{Filbet2002}
F.~Filbet and L.~Pareschi.
\newblock A numerical method for the accurate solution of the
  {F}okker-{P}lanck-{L}andau equation in the nonhomogeneous case.
\newblock {\em J. Comput. Phys.}, 179(1):1--26, 2002.

\bibitem{Filbet2016}
F.~Filbet and L.~M. Rodrigues.
\newblock Asymptotically {S}table {P}article-{I}n-{C}ell {M}ethods for the
  {V}lasov--{P}oisson {S}ystem with a {S}trong {E}xternal {M}agnetic {F}ield.
\newblock {\em SIAM J. Numer. Anal.}, 54(2):1120--1146, 2016.

\bibitem{filbet2003}
F.~Filbet and E.~Sonnendr{\"u}cker.
\newblock Comparison of {E}ulerian {V}lasov solvers.
\newblock {\em Computer Physics Communications}, 150(3):247--266, 2003.

\bibitem{Filbet2001}
F.~Filbet, E.~Sonnendr{\"u}cker, and P.~Bertrand.
\newblock Conservative {N}umerical {S}chemes for the {V}lasov {E}quation.
\newblock {\em J. Comput. Phys.}, 172(1):166--187, 2001.

\bibitem{HE81}
R.~W. Hockney and J.~W. Eastwood.
\newblock {\em Computer {S}imulation {U}sing {P}articles}.
\newblock McGraw Hill International Book Co., 1981.

\bibitem{Jingwei2018}
J.~Hu, S.~Jin, and R.~Shu.
\newblock A stochastic {G}alerkin method for the {F}okker-{P}lanck-{L}andau
  equation with random uncertainties.
\newblock In {\em Theory, numerics and {A}pplications of {H}yperbolic
  {P}roblems. {II}umerics and applications of hyperbolic problems. {II}},
  volume 237 of {\em Springer Proc. Math. Stat.}, pages 1--19. Springer, Cham,
  2018.

\bibitem{Jingwei2021}
J.~Hu, L.~Pareschi, and Y.~Wang.
\newblock Uncertainty quantification for the {BGK} model of the {B}oltzmann
  equation using multilevel variance reduced {M}onte {C}arlo methods.
\newblock {\em SIAM/ASA J. Uncertain. Quantif.}, 9(2):650--680, 2021.

\bibitem{landau1965}
L.~D. Landau.
\newblock On the vibrations of the electronic plasma.
\newblock {\em J. Phys.}, 10(1):25--34, 1965.

\bibitem{liu2017}
C.~Liu and K.~Xu.
\newblock A unified gas kinetic scheme for continuum and rarefied flows {V}:
  {M}ultiscale and multi-component plasma transport.
\newblock {\em Commun. Comput. Phys.}, 22(5):1175--1223, 2017.

\bibitem{McKinstrie1999}
C.~J. McKinstrie, R.~E. Giacone, and E.~A. Startsev.
\newblock Accurate formulas for the {L}andau damping rates of electrostatic
  waves.
\newblock {\em Phys. Plasmas}, 6:463--466, 1999.

\bibitem{medaglia2022}
A.~Medaglia, A.~Tosin, and M.~Zanella.
\newblock Monte {C}arlo stochastic {G}alerkin methods for non-{M}axwellian
  kinetic models of multiagent systems with uncertainties.
\newblock {\em Partial Differ. Equ. Appl.}, 3:51, 2022.

\bibitem{Pareschi2021}
L.~Pareschi.
\newblock An introduction to uncertainty quantification for kinetic equations
  and related problems.
\newblock In {\em Trails in {K}inetic {T}heory: {F}oundational {A}spects and
  {N}umerical {M}ethods}, volume~25 of {\em SEMA SIMAI Springer Ser.}, pages
  141--181. Springer, Cham, 2021.

\bibitem{Pareschi2001}
L.~Pareschi and G.~Russo.
\newblock An introduction to {M}onte {C}arlo method for the {B}oltzmann
  equation.
\newblock {\em ESAIM: Proc.}, 10:35--75, 2001.

\bibitem{pareschi2001_timerelaxed}
L.~Pareschi and G.~Russo.
\newblock Time relaxed {M}onte {C}arlo methods for the {B}oltzmann equation.
\newblock {\em SIAM J. Sci. Comput.}, 23(4):1253--1273, 2001.

\bibitem{Pareschi2000}
L.~Pareschi, G.~Russo, and G.~Toscani.
\newblock Fast spectral methods for the {F}okker-{P}lanck-{L}andau collision
  operator.
\newblock {\em J. Comput. Phys.}, 165(1):216--236, 2000.

\bibitem{Pareschi2013}
L.~Pareschi and G.~Toscani.
\newblock {\em Interacting {M}ultiagent {S}ystem: {K}inetic equations and
  {M}onte {Carlo} methods}.
\newblock Oxford University Press, 2013.

\bibitem{pareschi2005}
L.~Pareschi and S.~Trazzi.
\newblock Numerical solution of the {B}oltzmann equation by time relaxed
  {M}onte {C}arlo ({TRMC}) methods.
\newblock {\em Int. J. Numer. Methods Fluids}, 48(9):947--983, 2005.

\bibitem{Pareschi2020}
L.~Pareschi and M.~Zanella.
\newblock Monte {C}arlo stochastic {G}alerkin methods for the {B}oltzmann
  equation with uncertainties: {S}pace-homogeneous case.
\newblock {\em J. Comput. Phys.}, 423:1098--1120, 2020.

\bibitem{Poelle2019}
G.~Po\"{e}tte.
\newblock A g{PC}-intrusive {M}onte-{C}arlo scheme for the resolution of the
  uncertain linear {B}oltzmann equation.
\newblock {\em J. Comput. Phys.}, 385:135--162, 2019.

\bibitem{rossmanith2011}
J.~A. Rossmanith and D.~C. Seal.
\newblock A positivity-preserving high-order semi-lagrangian discontinuous
  {G}alerkin scheme for the {V}lasov--{P}oisson equations.
\newblock {\em J. Comput. Phys.}, 230(16):6203--6232, 2011.

\bibitem{shu2009}
C.-W. Shu.
\newblock High order weighted essentially nonoscillatory schemes for convection
  dominated problems.
\newblock {\em SIAM review}, 51(1):82--126, 2009.

\bibitem{Jin2019b}
R.~Shu and S.~Jin.
\newblock A study of {L}andau damping with random initial inputs.
\newblock {\em J. Differential Equations}, 266(4):1922--1945, 2019.

\bibitem{sod1978}
G.~A. Sod.
\newblock A survey of several finite difference methods for systems of
  nonlinear hyperbolic conservation laws.
\newblock {\em J. Comput. Phys.}, 27(1):1--31, 1978.

\bibitem{Sonnendrucker2013}
E.~Sonnendr{\"u}cker.
\newblock {\em Numerical methods for the {V}lasov equations}.
\newblock Lecture notes, 2013.

\bibitem{strang1968}
G.~Strang.
\newblock On the construction and comparison of difference schemes.
\newblock {\em SIAM J. Numer. Anal.}, 5(3):506--517, 1968.

\bibitem{Villani2013}
C.~Villani.
\newblock Landau damping.
\newblock In {\em Numerical models for fusion}, volume 39/40 of {\em Panor.
  Synth\`eses}, pages 237--326. Soc. Math. France, Paris, 2013.

\bibitem{Xiao2021}
T.~Xiao and M.~Frank.
\newblock A stochastic kinetic scheme for multi-scale plasma transport with
  uncertainty quantification.
\newblock {\em J. Comput. Phys.}, 432:110--139, 2021.

\bibitem{xiu2010}
D.~Xiu.
\newblock {\em Numerical {M}ethods for {S}tochastic {C}omputations}.
\newblock Princeton University Press, 2010.

\bibitem{xiu2002}
D.~Xiu and G.~E. Karniadakis.
\newblock The {W}iener--{A}skey polynomial chaos for stochastic differential
  equations.
\newblock {\em SIAM J. Sci. Comput.}, 24(2):619--644, 2002.

\bibitem{Zhang2017}
C.~Zhang and I.~M. Gamba.
\newblock A conservative scheme for {V}lasov {P}oisson {L}andau modeling
  collisional plasmas.
\newblock {\em J. Comput. Phys.}, 340:470--497, 2017.

\end{thebibliography}

\end{document}